\def\UseSection{
      \numberwithin{equation}{section}
	\theoremstyle{plain}
      \newtheorem{theorem}    {Theorem}[section]
      \DefineTheorems 
}

\def\DefineTheorems{
	
	\newtheorem{lemma}      [theorem] {Lemma}
	
	\newtheorem{prop}       [theorem] {Proposition}
	
	\newtheorem{cor}        [theorem] {Corollary}

	\theoremstyle{definition}
	\newtheorem{defn}       [theorem] {Definition}

	\theoremstyle{definition}

}

\newcommand{\bt}   {\begin{theorem}}
\newcommand{\et}   {\end  {theorem}}
\newcommand{\bl}   {\begin{lemma}}
\newcommand{\el}   {\end  {lemma}}
\newcommand{\bp}   {\begin{prop}}
\newcommand{\ep}   {\end  {prop}}
\newcommand{\bc}   {\begin{cor}}
\newcommand{\ec}   {\end  {cor}}
\newcommand{\bd}   {\begin{defn}}
\newcommand{\ed}   {\end  {defn}}

\newcommand{\ba}   {\begin{array}}
\newcommand{\ea}   {\end  {array}}
\newcommand{\be}   {\begin{enumerate}}
\newcommand{\ee}   {\end  {enumerate}}
\newcommand{\bi}   {\begin{itemize}}
\newcommand{\ei}   {\end  {itemize}}

\def\eq#1\en{\begin{equation}#1\end{equation}}  
\def\eqsplit#1\ensplit{
	\begin{equation}\begin{split}#1\end{split}\end{equation}
	}
\def\eqalign#1\enalign{
	\begin{align}#1\end{align}
	}
\def\eqmul#1\enmul{
	\begin{multline}#1\end{multline}
	}
\newcommand{\eqarrstar} {\begin{eqnarray*}} 
\newcommand{\enarrstar} {\end{eqnarray*}} 
\newcommand{\eqarray}   {\begin{eqnarray}} 
\newcommand{\enarray}   {\end{eqnarray}}

\newcommand{\lbeq}[1]  {\label{e:#1}}
\newcommand{\refeq}[1] {\eqref{e:#1}}    

\newcommand{\smallW}{\scriptstyle \leftarrow}

\newcommand{\ssmallW}{\scriptscriptstyle \leftarrow}
\newcommand{\smallN}{\scriptstyle \uparrow}
\newcommand{\smallS}{\scriptstyle \downarrow}

\newcommand{\smallnw}{\scriptscriptstyle \nwarrow}
\newcommand{\smallsw}{\scriptscriptstyle \swarrow}

\newcommand{\OTSP}{\begin{picture}(,)
\put(6.2,0){$\smallN$}
\put(.3,-3.3){$\smallW$}
\put(1.7,0){$\smallnw$}
\end{picture}\hspace{.4cm}
}
\newcommand{\smallFSOSP}{\begin{picture}(,)
\put(5.5,2.7){$\smallN$}
\put(5.5,-2.3){$\smallS$}
\put(.9,3.1){$\smallnw$}
\put(.9,-2){$\smallsw$}
\put(0.5,.5){$\ssmallW$}
\end{picture}\hspace{.35cm}
}

\makeatletter
\newcommand{\labelcounter}[2]{{%
	\stepcounter{#1}
	\protected@write\@auxout{}%
	{\string\newlabel{#2}{{\csname the#1\endcsname}{\thepage}}}%
	{\ref{#2}}
	}}
\makeatother
\newcommand{\sss}   { \scriptscriptstyle } 

\newcommand{\Nbold} {{\mathbb N}}

\newcommand{\Rbold} {{\mathbb R}}

\newcommand{\Zbold} {{\mathbb Z}}

\newcommand{\spose}[1] {{\hbox to 0pt{#1\hss}} }
\newcommand{\ltapprox} {\mathrel{\spose{\lower 3pt\hbox{$\mathchar"218$}}
\raise 2.0pt\hbox{$\mathchar"13C$}}}
\newcommand{\gtapprox} {\mathrel{\spose{\lower 3pt\hbox{$\mathchar"218$}}
\raise 2.0pt\hbox{$\mathchar"13E$}}}

\documentclass[12pt]{amsart}
\usepackage{amsmath,amssymb,amsfonts,amsthm,bbm}
\usepackage[dvips]{graphicx}
\newtheorem{THM}{Theorem}[section]

\newtheorem{COR}[THM]{Corollary}
\newtheorem{EXA}[THM]{Example}

\newtheorem{LEM}[THM]{Lemma}

\newtheorem{DEF}[THM]{Definition}

\UseSection  
\newcommand{\hlf}{\frac{1}{2}}
\newcommand{\ra}{\rightarrow}
\newcommand{\lra}{\leftrightarrow}

\newcommand{\la}{\leftarrow}

\renewcommand{\to}      {\rightarrow}

\newcounter{countC}  
\newcounter{countR}  
\newcommand{\R}{\Rbold}
\newcommand{\re}{\mathbb{R}}
\newcommand{\Z}{\Zbold}
\newcommand{\N}{\Nbold}

\newcommand{\mc}[1]{\mathcal{#1}}
\newcommand{\mP}{\mathbb{P}}

\newcommand{\mE}{\mathbb{E}}

\newcommand{\UD}{\updownarrow}
\newcommand{\LR}{\leftrightarrow}
\newcommand{\WE}{\LR}
\newcommand{\NS}{\UD}

\newcommand{\NE}{\begin{picture}(,)
\put(2,-5){$\rightarrow$}
\put(0,.5){$\uparrow$}
\end{picture}\hspace{.5cm}
}

\newcommand{\NEalpha}{\begin{picture}(,)
\put(2,-5){$\rightarrow$}
\put(4,-8){$\scriptstyle\alpha$}
\put(0,.5){$\uparrow$}
\end{picture}\hspace{.5cm}
}
\newcommand{\SE}{\begin{picture}(,)
\put(1.5,4.8){$\rightarrow$}
\put(-.5,-.5){$\downarrow$}
\end{picture}\hspace{.5cm}
}

\newcommand{\SEalpha}{\begin{picture}(,)
\put(2,5){$\rightarrow$}
\put(4,12){$\scriptstyle\alpha$}
\put(0,.3){$\downarrow$}
\end{picture}\hspace{.5cm}
}
\newcommand{\SW}{\begin{picture}(,)
\put(0,4.8){$\leftarrow$}
\put(8.2,-0.5){$\downarrow$}
\end{picture}\hspace{.5cm}
}

\newcommand{\NW}{\begin{picture}(,)
\put(0,-5){$\leftarrow$}
\put(8,0){$\uparrow$}
\end{picture}\hspace{.5cm}
}

\newcommand{\EWbeta}{
\underset{\hspace{3.5mm}\beta}{\longleftrightarrow}
}
\newcommand{\lowEWbeta}{\begin{picture}(,)
\put(.1,-5){$\longleftrightarrow$}
\put(13,-11){$\scriptstyle\beta$}
\end{picture}\hspace{.5cm}
}
\newcommand{\highEWbeta}{\begin{picture}(,)
\put(.1,5){$\longleftrightarrow$}
\put(13,0){$\scriptstyle\beta$}
\end{picture}\hspace{.5cm}
}
\newcommand{\NSW}{\begin{picture}(,)
\put(7.5,0){$\updownarrow$}
\put(0,0){$\leftarrow$}
\end{picture}\hspace{.5cm}
}
\newcommand{\NSE}{\begin{picture}(,)
\put(0,0){$\updownarrow$}
\put(2.2,0){$\rightarrow$}
\end{picture}\hspace{.5cm}
}
\newcommand{\SWE}{\begin{picture}(,)
\put(0,5){$\leftarrow$}
\put(5,5){$\rightarrow$}
\put(5.5,-0.5){$\downarrow$}
\end{picture}\hspace{.6cm}
}

\newcommand{\NWE}{\begin{picture}(,)
\put(0,-5){$\leftarrow$}
\put(5,-5){$\rightarrow$}
\put(5.5,0.5){$\uparrow$}
\end{picture}\hspace{.6cm}
}

\newcommand{\NSEWalt}{\begin{picture}(,)
\put(8,0){$\updownarrow$}
\put(0,0){$\longleftrightarrow$}
\end{picture}\hspace{.5cm}
}

\newcommand{\blank}[1]{}

\newcommand{\Qed}{\qed \medskip}

\usepackage[usenames]{color}
\newcommand{\red}{\color{red}}

\newcommand{\vep}{\varepsilon}

\newcommand{\ccup}[2]{\overset{\sss #2}{\underset{\sss #1}{\cup}}}
\newcommand{\ccap}[2]{\overset{\sss #2}{\underset{\sss #1}{\cap}}}

\begin{document}

\title[RWDRE]{Random walks in degenerate random environments.}

\author[Holmes]{Mark Holmes}
\address{Department of Statistics, University of Auckland}
\email{mholmes@stat.auckland.ac.nz}
\author[Salisbury]{Thomas S. Salisbury} 
\address{Department of Mathematics and Statistics, York University}
\email{salt@yorku.ca}


\maketitle

\begin{abstract}
We study the asymptotic behaviour of random walks in i.i.d.~random environments on $\Z^d$.  The environments need not be elliptic, so some steps may not be available to the random walker.  We prove a  monotonicity result for the velocity (when it exists) for any 2-valued environment, and show that this does not hold for 3-valued environments without additional assumptions.    We give a proof of directional transience and the existence of positive speeds under strong, but non-trivial conditions on the distribution of the environment.  
Our results include generalisations (to the non-elliptic setting) of 0-1 laws for directional transience, and in 2-dimensions the existence of a deterministic limiting velocity.  
\end{abstract}

\section{Introduction}
We will study simple random walks in random environments (RWRE) that are degenerate, in the sense that no ellipticity condition is assumed.  Our main results can be illustrated via the following example.  

\begin{figure}
\begin{center}
\includegraphics[scale=.45]{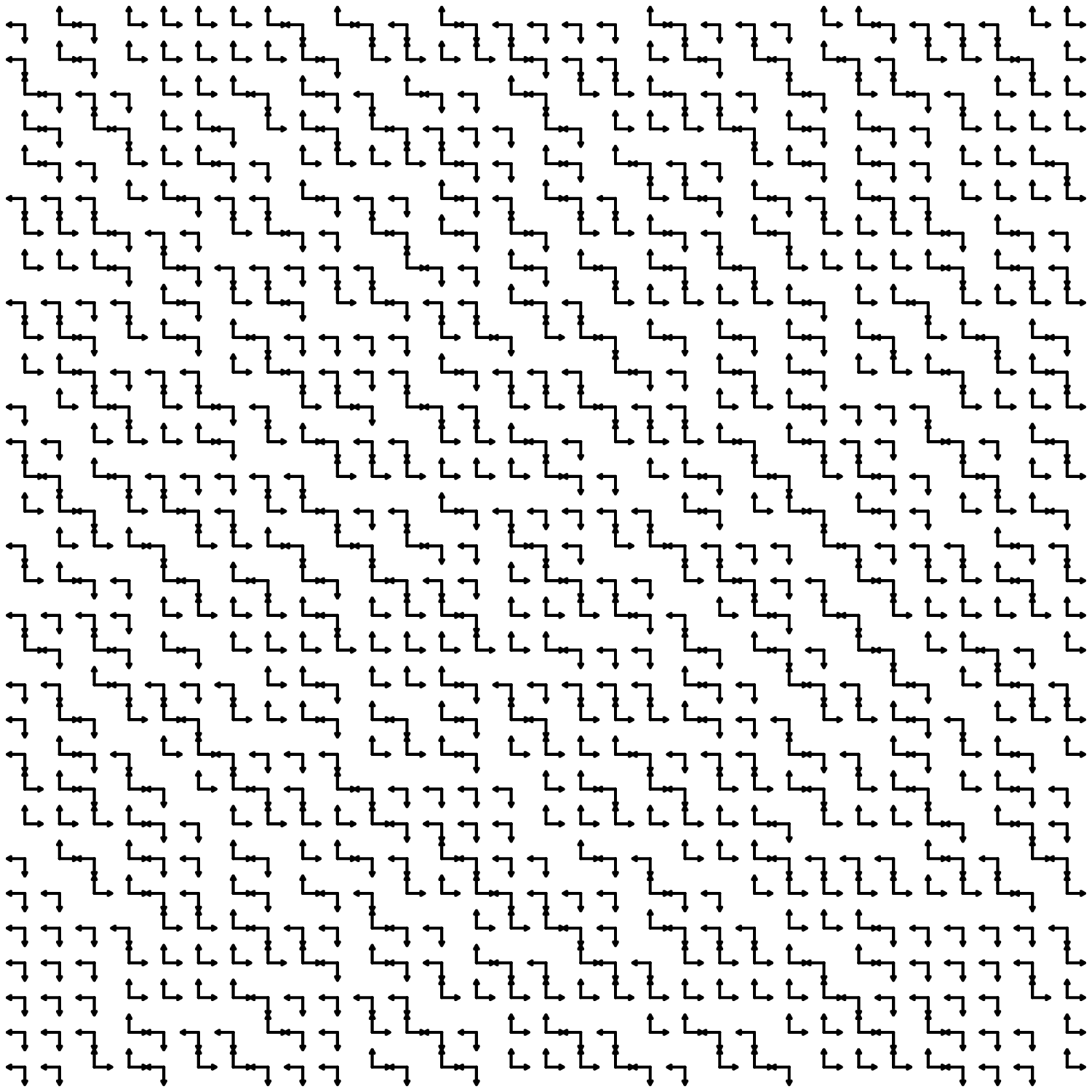} 
\includegraphics[scale=.45]{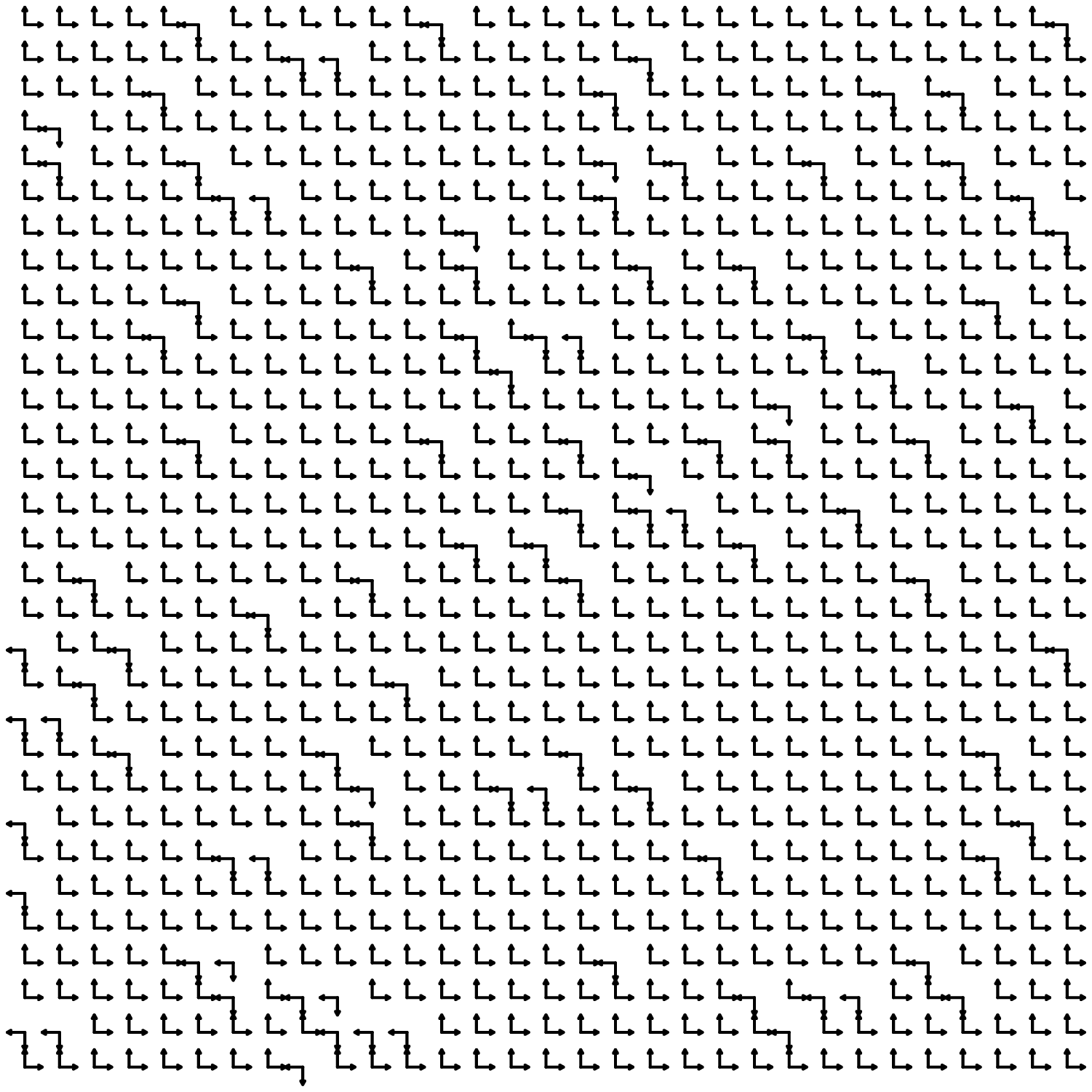} 
\end{center}
\caption{Finite regions of the random environment in Example~\ref{exa:NE_SW} for $p=.5$ and $p=.9$ respectively.}
\label{fig:orthant_env}
\end{figure}

\begin{EXA} $(\NE\SW)$:
\label{exa:NE_SW}
Perform site percolation with parameter $p$ on the lattice $\Z^2$.  From each occupied vertex $x=(x^{[1]},x^{[2]})$, insert two directed edges, one pointing up $\uparrow$ and one pointing right $\rightarrow$.  If $x$ is not occupied, insert directed edges pointing down $\downarrow$ and left $\leftarrow$ (see Figure \ref{fig:orthant_env}).  Now start a random walk at the origin $o$ that evolves by choosing uniformly from available arrows at its current location.
\end{EXA}

We conjecture that the random walk in this example has a speed in direction $\nearrow$ that is {\em strictly} monotone in $p\in [0,1]$, and that there is an infinite recurrent set of sites when $p=\hlf$.  We do not know how to prove these conjectures, however when the main results of this paper are applied to this model we will have established that:
\begin{itemize}
\item[(i)] for each $p\in [0,1]$, the walk visits infinitely many sites almost surely,
\item[(ii)] for $p$ sufficiently large, the random walk is transient in direction $\nearrow$,
\item[(iii)] for each $p$ the speed exists almost surely in the direction of $\ell$ for each $\ell \in \re^2\setminus o$,
\item[(iv)] the speed of the walk in direction $\nearrow$ is monotone increasing in $p$.
\end{itemize}
We will in fact show that (i) and (ii) follow from certain connectivity properties of a random directed graph determined by the environment, and studied in an earlier paper \cite{HS_DRE1}.  Result (iii) depends on an extension of a result of Zerner and Merkl (to the non-elliptic setting) that is valid for i.i.d.~RWRE models in 2-dimensions.  We will prove a version of (iv) that is valid for all i.i.d.~models where the local  environment has only 2 possible values (and show that when 3 values are allowed, monotonicity may fail). 

The above example is one of many interesting $2$-valued examples in which the random walker chooses uniformly from available steps at each site.  Our initial interest was in these models, so most of our examples will be of this kind, but the proofs apply more generally.

This paper is organised as follows.  In Section \ref{sec:model} we define random walks in random environments, and the directed graph associated to an environment, and state our main results.  In Section \ref{sec:RWREconnections} we recall from \cite{HS_DRE1,HS_DRE2} notions and results about connectivity in such random directed graphs, and 
examine RWRE results that can be inferred directly from the connectivity properties of these graphs.  In Section \ref{sec:orthog01} we adapt established techniques from the elliptic setting to prove 0-1 laws for directional transience and recurrence in our setting.  In Section \ref{sec:coupling} we use coupling methods to prove transience, ballisticity and monotonicity of speeds for certain models.  
Finally in Section \ref{sec:renewal} we give explicit speed formulae for some 2-valued 2-dimensional models with a simple renewal structure.


\subsection{The model}
\label{sec:model}

For fixed $d\ge 2$ let $\mc{E}=\{\pm e_i: i=1,\dots,d\}$ be the set of unit vectors in $\Z^d$.  Let $\mc{P}=M_1(\mc{E})$ denote the set of probability measures on $\mc{E}$, and let $\mu$ be a probability measure on $\mc{P}$.  If $\gamma\in \mc{P}$ we will abuse notation and write $\mu(\gamma)$ for $\mu(\{\gamma\})$.  Let $\Omega=\mc{P}^{\Z^d}$ be equipped with the product measure $\nu=\mu^{\otimes \Z^d}$ (and the corresponding product $\sigma$-algebra).  An environment $\omega=(\omega_x)_{x\in \Z^d}$ is an element of $\Omega$.  We write $\omega_x(e)$ for $\omega_x(\{e\})$. Note that $(\omega_x)_{x\in \Z^d}$ are i.i.d.~with law $\mu$ under $\nu$.

The random walk in environment $\omega$ is a time-homogeneous (but not necessarily irreducible) Markov chain with transition probabilities from $x$ to $x+e$ defined by 
\begin{equation}
\label{eq:trans_prob}
p_{\omega}(x,x+e)=
\omega_x(e).
\end{equation}
Given an environment $\omega$, we let $\mP_{\omega}$ denote the law of this random walk $X_n$, starting at the origin.  Let $P$ denote the law of the annealed/averaged random walk, i.e.~$P(\cdot, \star):=\int_{\star}\mP_{\omega}(\cdot)d\nu$.
Since $P(A)=E_{\nu}[\mP_{\omega}(A)]$ and $0\le \mP_{\omega}(A)\le 1$, $P(A)=1$ if and only if $\mP_{\omega}(A)=1$ for $\nu$-almost every $\omega$.  Similarly $P(A)=0$ if and only if $\mP_{\omega}(A)=0$ for $\nu$-almost every $\omega$.  
If we start the RWRE at $x\in\Z^d$ instead, we write $P_x$ for the corresponding probability, so $P=P_o$.

We associate to each environment $\omega$ a directed graph $\mc{G}(\omega)$ (with vertex set $\Z^d$) as follows.  For each $x\in \Z^d$, the directed edge $(x,x+u)$ is in $\mc{G}_x$ if and only if $\omega_x(u)>0$, and the edge set of $\mc{G}(\omega)$ is $\cup_{x\in \Z^d} \mc{G}_x(\omega)$.  For convenience we will also write $\mc{G}=(\mc{G}_x)_{x\in \Z^d}$.    
Note that under $\nu$, $(\mc{G}_x)_{x\in \Z^d}$ are i.i.d.~subsets of $\mc{E}$.  The graph $\mc{G}(\omega)$ is equivalent to the entire graph $\Z^d$ (with directed edges), precisely when the environment is {\em elliptic}, i.e.~$\nu(\omega_x(u)>0)=1$ for each $u \in \mc{E}, x\in \Z^d$.  Much of the current literature assumes either the latter condition, or the stronger property of {\em uniform ellipticity}, i.e.~that $\exists \epsilon>0$ such that $\nu(\omega_x(u)>\epsilon)=1$ for each $u \in \mc{E}, x\in \Z^d$.

On the other hand, given a directed graph $\mc{G}=(\mc{G}_x)_{x\in \Z^d}$ (with vertex set $\Z^d$, and such that $\mc{G}_x\ne \varnothing$ for each $x$), we can define a {\em uniform} random environment $\omega=(\omega_x(\mc{G}_x))_{x\in \Z^d}$. Let $|A|$ denote the cardinality of $A$, and set
\[\omega_x(e)=\begin{cases}
|\mc{G}_x|^{-1}, & \text{ if }e \in \mc{G}_x\\
0, & \text{otherwise}.
\end{cases}\]
The corresponding RWRE then moves by choosing uniformly from available steps at its current location.  This gives us a way of constructing rather nice and natural examples of random walks in non-elliptic random environments:  first generate a random directed graph $\mc{G}=(\mc{G}_x)_{x\in \Z^d}$ where $\mc{G}_x$ are i.i.d., then run a random walk on the resulting random graph (choosing uniformly from available steps).  
This natural class of RWRE will henceforth be referred to as {\em uniform RWRE}.  Note that we have chosen above to forbid $\mc{G}_x=\varnothing$.  In the setting of uniform RWRE it would be reasonable to instead allow $\mc{G}_x=\varnothing$ and to define $\omega_x(o)=1$ in this case, with the walker getting absorbed at $x$. However (see Lemma \ref{lem:stuck} below) if this happens with positive probability then the random walker gets stuck on a finite set of vertices almost surely, so in terms of the random walk behaviour, we lose no interesting cases by prohibiting $\mc{G}_x=\varnothing$.  More generally, we have the following explicit criterion for whether a RWRE gets stuck on a finite set of sites, which depends only on $\omega$ via the connectivity of $\mc{G}(\omega)$.  Here $V$ is an orthogonal set if $u\cdot v=0$ for every $u,v\in V$.

\begin{THM} 
\label{thm:stuck}
If there exists a nonempty orthogonal set $V\subset \mc{E}$ such that $\mu(\mc{G}_o\cap V\neq\varnothing)=1$ then the random walk visits infinitely many sites, $P$-almost surely.  Otherwise the random walk visits only finitely many sites $P$-almost surely.
\end{THM}

Fix $\ell\in \re^d\setminus o$.  Let $A_+$ and $A_-$ denote the events that $X_n\cdot \ell\to\infty$ and $X_n\cdot \ell\to-\infty$ respectively.  Clearly when the random walk gets stuck on a finite set of sites, the walker is not directionally transient in any direction and the speed of the walk is zero.  By examining the connectivity structure of $\mc{G}$ we are able to prove the following generalisation (to the non-elliptic setting) of a 0-1 law first proved in the uniformly elliptic setting by Kalikow \cite{K81}. 

\begin{THM}
\label{thm:A+A-}
For random walks in i.i.d.~random environments, $P(A_+\cup A_-)\in \{0,1\}$.
\end{THM}
This allows us to extend the following, obtained by Sznitman and Zerner \cite{SZ99}, Zerner \cite{Zern02,Zern07} and Zerner and Merkl \cite{MZ01} in the elliptic setting, to random walks in i.i.d.~random environments with no ellipticity assumption.
\begin{THM}
\label{thm:dtransv}
There exist deterministic $v_{+}(\ell), v_{-}(\ell)$ such that 
\[\lim_{n\ra \infty}\frac{X_n\cdot \ell}{n}=v_{+}(\ell)\mathbbm{1}_{A_+}+v_{-}(\ell)\mathbbm{1}_{A_{-}}, \quad P-\text{a.s.}\]
\end{THM}
By applying Theorem \ref{thm:dtransv} to each of the standard basis vectors $\{e_j\}_{j=1}^d$, we have that $X_n/n$ has a limit $P-$almost surely.  In principle this limit could be random, taking at most 2 possible values.  
\begin{THM}
\label{thm:zern_merk}
When $d=2$, $P(A_{\ell})\in \{0,1\}$.
\end{THM}
These two theorems imply (as claimed above) that a deterministic velocity $v$ always exists in the 2-dimensional setting (see Corollary \ref{cor:speedsin2d}).  
We believe that this holds for all $d$.  This is simple in settings where renewals occur due to a forbidden direction, such as if $\mu(-e_1\in \mc{G}_o)=0$ but $\mu(e_1\in \mc{G}_o)>0$, and indeed much more can be said in such cases, see e.g.~\cite{RS06} and Section \ref{sec:renewal} and Table \ref{tab:walks}.  Of course, existence of speeds does not imply transience unless one can prove that $v\ne o$. When there is sometimes a drift in direction $u$ but never a drift in direction $-u$ the walk should be almost surely transient in direction $u$.
Some results of this kind are known in the uniformly elliptic setting, and we hope to address these issues in a subsequent paper without the assumption of ellipticity.  Instead, in this paper we give a relatively simple proof under the strong assumption that with sufficiently large probability we have a sufficiently large drift at the origin (see Theorem \ref{thm:excouple}), relying on results from \cite{HS_comb, Zern05, BS08}.  



\begin{DEF}
\label{def:2-val}
An environment is {\em $2$-valued} when there exist distinct $\gamma^{\sss 1},\gamma^{\sss 2}\in \mc{P}$ and $p \in (0,1)$ such that $\mu(\gamma^{\sss 1})=p$, $\mu(\gamma^{\sss 2})=1-p$.  A graph is {\em $2$-valued} when there exist distinct $E^{\sss 1},E^{\sss 2}\subset \mc{E}$ and $p\in(0,1)$ such that $\mu(\mc{G}_o=E^{\sss 1})=p$ and $\mu(\mc{G}_o=E^{\sss 2})=1-p$.  We write $(\gamma^{\sss 1},\gamma^{\sss 2})$ or $(E^{\sss 1},E^{\sss 2})$ to denote the family of 2-valued environments or graphs indexed by $p=\mu(\gamma^{\sss 1})$ or $p=\mu(E^{\sss 1})$.
\end{DEF}
Note that when an environment is 2-valued, the corresponding graph is at most 2-valued.  When a graph is 2-valued, the {\em uniform} environment corresponding to that graph is also 2-valued.  The uniform RWRE $(\NE,\SW)$ (see Example \ref{exa:NE_SW}) is an interesting 2-valued 2-dimensional example that has a natural generalisation to higher dimensions.  There are many other interesting 2-valued examples, such as the 2-dimensional uniform random environments $(\SWE,\uparrow)$, $(\NSEWalt \,\, ,\NE)$, $(\NS,\WE)$.  The last of these is a degenerate version of the ``good-node bad-node'' model of Lawler \cite{Law82}, which has been studied recently by Berger and Deuschel \cite{BD11}.  

The following theorem is one of the main results of this paper.
\begin{THM}
\label{thm:main1}
Fix any $2$-valued environment with $\mu(\gamma^{\sss 1})=p$, $\mu(\gamma^{\sss 2})=1-p$.  If for every $p$ there exists $v[p]$ such that $P(v[p]=\lim_{n\ra \infty}n^{-1}X_n)=1$, then each coordinate of $v[p]$ is monotone in $p$.
\end{THM}
In fact we'll prove a more comprehensive version of this in Section \ref{sub:coupling_mono}, via a simple coupling argument. See Theorem \ref{thm:other_mono}.  Note that this implies monotonicity of the speed as a function of $p$ for all 2-valued models in 2 dimensions.  We also show that the corresponding statement for $3$-valued models with probabilities $(p(1-q),(1-p)(1-q),q)$ (where $q$ is fixed) fails in general.

\section{Random walk properties obtained from the geometry of connected clusters}
\label{sec:RWREconnections}

\begin{DEF}
Given a directed graph $\mc{G}$:
\begin{itemize}
\item We say that $x$ is {\em connected} to $y$, and write $x\ra y$ if: there exists an $n\ge 0$ and a sequence $x=x_0\,x_1,\dots, x_n=y$ such that $x_{i+1}-x_{i}\in \mc{G}_{x_i}$ for $i=0,\dots,n-1$.  Let $\,\mc{C}_x=\{y\in \Z^d:x \ra y\}$, and $\mc{B}_y=\{x\in\Z^d:x\ra y\}$. 
\item We say that $x$ and $y$ are {\em mutually connected}, or that they {\em communicate}, and write $x\lra y$ if $x\ra y$ and $y \ra x$.  Let $\mc{M}_x=\{y \in \Z^d:x\lra y\}=\mc{B}_x\cap\mc{C}_x$. 
\item A cluster $\mc{M}$ is said to be {\em gigantic} in $\Z^d$ if every $\Z^d$-connected component of $\Z^d\setminus \mc{M}$ is finite.
\end{itemize}
For the $\mc{G}$ arising from the model of Section \ref{sec:model}, set
\[\theta_+=\nu(|\mc{C}_o|=\infty), \quad \theta_-=\nu(|\mc{B}_o|=\infty), \quad \text{and}\quad  \theta=\nu(|\mc{M}_o|=\infty).\]
\end{DEF}

In this section we present a number of results for RWRE, that depend only on the clusters $(\mc{C}_x)_{x\in \Z^d}$ and $(\mc{M}_x)_{x\in \Z^d}$ of the graph $\mc{G}(\omega)$ induced by the environment $\omega$.

Whether the RWRE $X$ gets stuck on a finite set of sites can be characterized completely in terms of the law of the connected cluster $\mc{C}_o$.  If $\mc{C}_o$ is almost surely infinite, then so is $\mc{C}_x$ for each $x$, so the random walk will eventually escape from any finite set of sites.  On the other hand if $\mc{C}_o$ is finite with positive probability, then we will see that there is some $\delta>0$ such that each time the walk reaches a new $\|X_n\|_{\infty}$ maximum it has probability at least $\delta$ of being at a site with finite $\mc{C}$, whence the walk will eventually get stuck.  These arguments are formalised in the following result.
\begin{LEM}
\label{lem:stuck}
Fix $d\ge 2$, and let $X_n$ be the random walk in i.i.d.~random environment.
\begin{itemize}
\item[(i)] $\theta_+=1\Rightarrow P(\sup_{n\ge 1}|X_n|<\infty)=0$ (i.e.~the RWRE visits infinitely many sites).
\item[(ii)] $\theta_+<1\Rightarrow P(\sup_{n\ge 1}|X_n|<\infty)=1$ (i.e.~the RWRE gets stuck on a finite set of points).
\end{itemize}
\end{LEM}
\proof To prove (i), note that $\nu$-a.s., $\mc{C}_x(\omega)$ is infinite for every $x$.  Now  fix such an $\omega$ and 
assume that $X_n$ visits only finitely many sites. Then it must visit some site $y$ infinitely often. Let $z\in\mc{C}_y$. There is an admissible path connecting $y$ to $z$ which has a fixed positive probability of being followed on any excursion of $X_n$ from $y$. Therefore 
$X_n$ will eventually visit $z$.  
Since $\mc{C}_y$ is infinite $\nu$-a.s., this contradicts the assumption. 

To prove (ii), suppose that 
$\theta_+<1$.  Define
$$
n_0=\inf\{n\ge 1:\exists F\subset\Z^d \text{ with }|F|=n \text{ and }\nu(\mc{C}_o=F)>0\}
$$
and choose $F$ satisfying $|F|=n_0$ and $\delta=\nu(\mc{C}_o=F)>0$. 
Note that if $\mc{M}_y=F$ for some $y\in F$, then $\mc{M}_{y'}=F$ for each $y'\in F$.  So by translation invariance of $\nu$, for each $y\in F$,
\begin{equation}
\label{fmla:translationofM}
\nu(\mc{M}_o=F)=\nu(\mc{M}_y=F)=\nu(\mc{M}_o=F-y).
\end{equation}
Furthermore, for $y \in F$, if $F=\mc{C}_y\supsetneq\mc{M}_y$, then there exists $y'\in \mc{C}_y\setminus \mc{M}_y$.  Since $y \notin \mc{C}_{y'}\subset \mc{C}_y$ the set $G=\mc{C}_{y'}\subsetneq F$ satisfies $|G|<n_0$, and 
$\nu(\mc{C}_o=G-y')=\nu(\mc{C}_{y'}=G)>0$.  This would contradict the definition of $n_0$. 
So in fact $\nu(\mc{M}_y=F)=\nu(\mc{C}_y=F)$ for each $y\in F$. Therefore by (\ref{fmla:translationofM}), 
$$
\nu(\mc{C}_o=F)=\nu(\mc{C}_o=F-y),
$$
for every $y\in F$. 

For each $x\in \Z^d$ we can find a $y\in F$ such that $x$ minimizes $\|z\|_\infty$ over $z\in F+x-y$ (just find a unit vector $e\in\Z^d$ such that $\|x+e\|_\infty>\|x\|_\infty$, and then choose $y$ so that the projection of $F-y$ in the direction of $e$ is $\ge 0$). Then 
\begin{equation}
\nu(\mc{C}_x=F+x-y)=\nu(\mc{C}_o=F-y)=\nu(\mc{C}_o=F)=\delta.
\label{fmla:yonedgeofF}
\end{equation}

For $k\ge 1$ define
\[
T_k:=\inf\{m\ge 1:\|X_m\|_{\infty}=k(n_0+1)\}.
\]
Let $\mc{F}_k$ reveal $X_m$ for $m< T_k$, and the environment $\omega_z$ for $\|z\|_\infty<k(n_0+1)$. 
Then 
$$
P(T_k=\infty\mid\mc{F}_{k-1})\ge P(|\mc{C}_{X_{T_{k-1}}}|\le n_0 \mid\mc{F}_{k-1})
\ge\sum_{y\in F}P(\mc{C}_{X_{T_{k-1}}}=X_{T_{k-1}}+F-y\mid\mc{F}_{k-1}).
$$
We have shown above that $y$ can be chosen so that $X_{T_{k-1}}+F-y$ is disjoint from $\{z:\|z\|_\infty<(k-1)(n_0+1)\}$ (the region whose environment is revealed by $\mc{F}_{k-1}$). Applying  (\ref{fmla:yonedgeofF}) to that $y$ we conclude that
$$
P(T_k<\infty\mid\mc{F}_{k-1})\le (1-\delta)\mathbbm{1}_{\{T_{k-1}<\infty\}}.
$$
Iterating $k$ times, $P(T_k<\infty)\le (1-\delta)^k$, and sending $k\to\infty$ we obtain 
\[0=P\left(\cap_{k=1}^{\infty}\{T_k<\infty\}\right)=1-\nu\big(\cup_{k=1}^{\infty}\{T_k=\infty\}\big),\]
which establishes the result.
\Qed

On the event that the walk gets stuck on a finite set of sites, the asymptotic velocity is trivially zero and the walk is not directionally transient in any direction, almost surely.  Hence, by Lemma \ref{lem:stuck}, Theorems \ref{thm:A+A-}, \ref{thm:dtransv} and \ref{thm:zern_merk} hold trivially when 
$\theta_+<1$. 
Our principal interest will therefore be in situations where the following condition holds: 
\begin{equation}
\label{standinghypothesis}
\theta_+=1.
\end{equation}
Note that our general hypotheses rule out the possibility that $\mc{C}_x=\{x\}$. As remarked above, at the cost of more cumbersome notation, we could have included the possibility of $\mc{G}_x=\varnothing$ in our models. In this case also 
$\theta_+<1$ 
and the proof of Lemma \ref{lem:stuck}(ii) remains valid.
\medskip

Together with Lemma \ref{lem:stuck}, the following simple criterion from \cite{HS_DRE1} 
proves Theorem \ref{thm:stuck}.

\begin{LEM} 
$\theta_+=1$ 
if and only if there exists an orthogonal set $V$ of unit vectors such that $\mu(\mc{G}_o\cap V\neq\varnothing)=1$.
\label{lem:criteriontobeinfinite}
\end{LEM}
It follows immediately that for the uniform RWRE (and indeed RWRE for any environment giving rise to such graphs) $(\NE,\SW)$ of Example \ref{exa:NE_SW}, and similarly for the models $(\NS,\WE)$, $(\SWE,\uparrow)$, the walk visits infinitely many sites almost surely, by choosing e.g.~$V=\{\uparrow,\leftarrow\}$ in each case.  On the other hand if for example $\{\NE,\NW, \SW, \SE\}$ is a subset of the possible values for $\mc{G}_o$ under $\mu$ then the walk has finite range. 

From \cite{HS_DRE1} (see also \cite{HS_DRE2} for improvements to some of these values), the following Lemma immediately implies that the RWRE is transient in the following situations: when $\mc{G}_o(\omega)\in \{\NE,\SW\}$  almost surely (as in Example~\ref{exa:NE_SW}), with $\mu(\mc{G}_o=\NE)>.83270$; 
when $\mc{G}_o(\omega)\in \{\SWE,\uparrow\}$  almost surely with $\mu(\mc{G}_o=\uparrow)>.83270$.

\begin{LEM}
\label{lem:Mtrans}
For any environment $\omega$ such that $|\mc{C}_x|=\infty$ and $|\mc{M}_x|<\infty$ for every $x \in \Z^d$, the random walk in environment $\omega$ is transient $\mP_{\omega}$-almost surely.
\end{LEM}
\proof Starting from any $x$, we visit infinitely many sites. Thus the time $T_x$ taken for the walk to exit $\mc{M}_x$ is finite, and any vertex reached thereafter is necessarily in $\mc{C}_x\setminus \mc{M}_x$.  Hence the random walk never returns to $x$ after time $T_x$.\Qed\medskip

For some of the most interesting models, such as the uniform RWRE $(\NS,\WE)$, $\mc{C}_o$ is infinite almost surely while $\mc{M}_o$ can be finite or infinite. 
In particular, the Markov chain in a typical fixed environment $\omega$ the Markov chain is not irreducible.  
Therefore, when considering recurrence, one should ask if the origin is visited infinitely often, given that $\mc{M}_o$ is infinite and $\mc{M}_o=\mc{C}_o$.  

\begin{LEM}
\label{lem:Mrecurrent}
For any environment $\omega$ such that $|\mc{C}_y|=\infty$ for every $y$ and such that there is a unique infinite mutually connected cluster $\mc{M}_{\infty}$, we have that $\{x:X_n=x \text{ infinitely often}\}\in \{\varnothing, \mc{M}_{\infty}\}$, $\mP_{\omega}$-almost surely.
\end{LEM}
\proof Let $\mc{R}=\{x:X_n=x \text{ infinitely often}\}$ be the recurrent set.  If $x\in \mc{R}$, then $\mc{C}_x\subset \mc{R}$ by the argument for (i) of Lemma \ref{lem:stuck}. Clearly also $\mc{R}\subset \mc{M}_{x}$ (every site in $\mc{R}$ must be reachable from $x$ and vice versa), almost surely.  Since $\mc{M}_x\subset \mc{C}_x$ this implies that $\mc{C}_x=\mc{R}=\mc{M}_x$ and also that $|\mc{M}_x|=\infty$ (so $\mc{M}_x=\mc{M}_{\infty}$).
\Qed

The above argument also shows that if there are multiple distinct infinite $\mc{M}_x$ then $\mc{R}$ is either empty or is equal to exactly one such infinite $\mc{M}_y$.  For the uniform RWRE $(\NS,\WE)$, \cite{HS_DRE1} shows that a unique cluster $\mc{M}_\infty$ exists $\nu$-a.s., as in Lemma \ref{lem:Mrecurrent}, for every $p\in(0,1)$. Combining recent results of Berger and Deuschel \cite{BD11} together with the argument of Theorem 3.3.22 of \cite{Zeit04} (attributed there to Kesten) it follows that in this model, $\mc{R}=\mc{M}_\infty$, $P$-a.s..

Another simple result is the following, which by results of \cite{HS_DRE1} applies to the models $(\NE,\SW)$ and $(\SWE, \uparrow)$ for certain values of $p$.
\begin{LEM}
\label{lem:giantMrecurrent}
Let $\omega$ satisfy $|\mc{C}_y|=\infty$ for every $y$ and suppose there exists $x$ such that $\mc{M}=\mc{M}_x$ is gigantic.  Then $T_{\mc{M}}\equiv\inf\{n\ge 0:X_n\in \mc{M}\}<\infty$ and $X_n\in \mc{M}$ for all $n\ge T_{\mc{M}}$, $\mP_{\omega}$-almost surely.
\end{LEM}

We now turn to a result that allows us to improve on transience and prove {\em directional} transience in some situations.  Given a vector $v\in \re^d\setminus o$, and $N\ge 1$, define $A_N^-(v)=\{y:y\cdot v<-N\}$. 

\begin{THM}
\label{thm:direction_trans}
Assume that $\theta_+=1$, and suppose that there is a vector $v\in \re^d\setminus o$ such that
\begin{itemize}
\item[(a)] $\nu\big(\cup_{N\ge 1} \{\mc{C}_o \cap A_N^-(v)=\varnothing\}\big)=1$, and
\item[(b)] $\mu(\{\exists e\in \mc{G}_o:v\cdot e>0)>0$.
\end{itemize}
Then $P(\liminf X_n\cdot v=\infty)=1$, i.e.~the random walk is a.s.~transient in direction $v$.
\end{THM}
\proof   
\noindent Without loss of generality, $|v|=1$.  For $x\in \Z^d$ define $F_x(v):=\{y:(y-x)\cdot v<0\}$ and ${\bf C}_v:=\{x: \mc{C}_x\cap F_x(v)=\varnothing\}$.
We show that under the hypotheses of the theorem:
\begin{itemize}
\item[(i)] $\nu(o \in {\bf C}_v)=\varepsilon$ for some $\varepsilon>0$
\item[(ii)] $P\big(\sup_{n}X_n\cdot v=\infty\big)=1$
\item[(iii)] for all $M<\infty$, $P\big(\liminf_{n\ra \infty} X_n\cdot v\ge M\big)=1$.
\end{itemize}
Note that the desired result clearly follows from (iii).

To prove (i), note that the result is trivial if for every $e \in \mc{E}$ such that $e\cdot v<0$ we have $\mu(e \in \mc{G}_o)=0$. 
So assume otherwise. Then there exist $e\in \mc{E}$, $\epsilon,\delta>0$ such that $e\cdot v\le -\delta$ and $\mu(\omega_o(e)>\epsilon)>\epsilon$.  Assume that (i) fails.  Then we can construct a sequence $o=x_0,x_1,\dots$ such that $x_i\in \mc{C}_{x_{i-1}}\cap F_{x_{i-1}}(v)$. Moreover we can identify $x_i$ without looking at the environment at any site $z$ such that $(z-x_i)\cdot v<0$.  Therefore any such site has probability at least $\epsilon^{2N}$ of being connected to a site that is at least distance $N\delta$ in the direction of $-v$ (from o).  We have infinitely many independent trials in which to observe such a connection. Thus almost surely we are connected at least distance $N\delta$ in the direction $-v$, for each $N$, contradicting (a).


To prove (ii), let $W_n=X_n\cdot v$ and let $U_v$ denote the set of unit vectors $e\in \Z^d$ such that $v\cdot e\ge 0$.  For any record levels $(z_e)_{e\in U_v}$ and any $N\ge 0$, consider
$$
\{w\in\Z^d: w\cdot v\ge -N\text{ and }w\cdot e\le z_e\,\forall e\in U_v\}.
$$
We claim this is finite. To see this, write $w=\sum w_iu_i$, where we select a basis $u_i$ from $U_v$. If $u_i\cdot v=0$, then $\pm u_i\in U_v$ and we have $z_{-u_i}\le w_i\le z_{u_i}$. If $u_i\cdot v>0$ then  $-N\le w\cdot v=\sum w_ju_j\cdot v\le \sum_{j\neq i} z_{u_j}u_j\cdot v+w_iu_i\cdot v$ so in fact 
$$
z_{u_i}\ge w_i\ge \frac{-N-\sum_{j\neq i} z_{u_j}u_j\cdot v}{u_i\cdot v}.
$$
In other words, there are only finitely many possibilities for each $w_i$.

By Lemma \ref{lem:stuck}, $X_n$ visits infinitely many sites.  
Since $\inf_n W_n>-\infty$ (which follows from the fact that $\nu(\cup_{N\ge 1} \{\mc{C}_o \cap A_N^-(v)=\varnothing\})=1$), 
 the above argument shows that the random times $T_1=\inf\{n\ge 1:\exists e\in U_v;X_n\cdot e>X_m\cdot e \text{ for all }m< n\}$ and 
 \[T_i=\inf\{n>T_{i-1}:\exists e \in U_v;X_n\cdot e>X_m\cdot e \text{ for all }m< n\}, \quad i=2,3,\dots\]
  are almost-surely finite, since eventually any set of record levels will be surpassed. 
  
Let  $V_x=\{y:(y-x)\cdot e\ge 0 \, \forall \, e \in U_v\}$. Then each time $T_i$ is the first hitting time of $V_{X_{T_i}}$, so the environment in $V_{X_{T_i}}$ is unexplored. Note that $(b)$ implies that there exists $\epsilon_v>0$ and $e \in U_v$ with $e\cdot v>0$ such that $\mu(\omega_o(e)>\epsilon_v)>\epsilon_v$. Set $U^0_v=\{e\in U_v:e\cdot v>0, \mu(\omega_o(e)>\epsilon_v)>\epsilon_v\}$, and $c_v=\{v\cdot e: e\in U^0_v\}>0$.  Then for each $k\in \mathbb{N}$, there is probability at least $\epsilon_v^{k}$ of there being an admissible path from $X_{X_{T_i}}$ consisting of $k$ arrows from $U^0_v$. Following this path keeps us in the previously unexplored region $V_{X_{T_i}}$.  So given this, there is probability at least $\epsilon_v^{k}$ that the next $k$ steps of $X_n$ will follow this path. In other words, for each $i,k\in \N$, there is probability at least $\epsilon_v^{2k}$ under $P$, independent of the history of the walk up to time $T_i$, that $\{(X_{T_i+k}-X_{T_i})\cdot v>c_vk\}$.  Therefore for each $k\in\N$ there will eventually be an $i$ such that $\{(X_{T_i+k}-X_{T_i})\cdot v>c_vk\}$. Recalling that $\inf W_n>-\infty$, we conclude that $\sup X_n\cdot v=\infty$.

To prove (iii), fix $M_1\ge 1$ and define $\mc{T}_1$ to be the first hitting time of $M_1$ by $W$, i.e.
\[\mc{T}_1=\inf\{n>0:W_n\ge M_1\},\]
which is $P$-almost surely finite by (ii). Given $\mc{T}_i<\infty$, 
let 
\[N_i=\inf\{k\ge 1: \text{ there exists a $\mc{G}$-admissible path (of length $k$) }  \vec{\eta}_k:X_{\mc{T}_i}\ra F_{X_{\mc{T}_i}}(v)\}.\] 
If $N_i<\infty$ set $M_{i+1}=M_i+N_i$ and 
\[\mc{T}_{i+1}=\inf\{n>\mc{T}_{i}:W_n\ge M_{i+1}\}.\]
Note that $\{N_i=\infty\}=\{X_{\mc{T}_i}\in {\bf C}_v\}$ and if $N_{i}=\infty$ then $W_n\ge M_i$ for all $n\ge \mc{T}_i$.  Moreover, to determine if $N_i\le m$, we need only look at the environment within distance $m$ of $X_{\mc{T}_i}$, so if $N_i<\infty$ then the walk visits an unexplored environment at time $\mc{T}_{i+1}$. In other words, the event that $X_{\mc{T}_i} \in {\bf C}_v$ depends only on the unexplored environment in $\Z^d\setminus F_{\mc{T}_i}(v)$, so that by (i), $\mc{I}=\inf\{i: X_{\mc{T}_i}\in {\bf C}_v\}$ has a geometric distribution with parameter $\varepsilon>0$. 

Thus $P$-almost surely there is a $\mc{T}_i<\infty$ such that $W_n\ge M_i\ge M_1$ for all $n\ge \mc{T}_i$.  Since $M_1$ was arbitrary, this proves (iii).
\Qed

\begin{COR}
\label{cor:Mdirect}
Assume that $\theta_+=1$.
For each $d\ge 2$ there exists $\epsilon_d$ such that the following holds:  
If there exists an orthogonal set $V$ of unit vectors such that $\mu(\mc{G}_o\subset V)>1-\epsilon_d$, then the random walk is transient in direction $v=\sum_{e\in V}e$.
\end{COR}
\proof The proof of \cite[Theorem 4.2]{HS_DRE1} verifies Theorem  \ref{thm:direction_trans} (a) for $v=\sum_{e\in V}e$, while condition (b) holds with $\epsilon_v=1-\epsilon_d$ by the assumption that $\mu(\mc{G}_o\subset V)>1-\epsilon_d$.\Qed


When applied to the models $(\NE \SW)$ and $(\SWE\uparrow)$ and using \cite[Corollary 4.3]{HS_DRE1}, Corollary \ref{cor:Mdirect} improves our transience results to transience in directions $(1,1)$ and $(0,1)$ when $\mu(\mc{G}_o=\NE)>0.83270$ and $\mu(\mc{G}_o=\uparrow)>0.83270$ respectively.  The former result is improved to $\mu(\mc{G}_o=\NE)>0.7491$ via the following Corollary.

\begin{COR}
\label{cor:orthant_D_trans}
Any 2-dimensional RWRE model with 2-valued graph $\mc{G}_o\in \{\NE, \SW\}$ is transient in direction $(1,1)$ when $\mu(\mc{G}_o=\NE)> p_c^{\OTSP}$, where $p_c^{\OTSP}$ is the critical occupation prob.~for oriented site percolation on the triangular lattice. 
\end{COR}
\proof As in \cite{HS_DRE2}, when $p=\mu(\mc{G}_o=\NE)>p_c^{\OTSP}$, $\mc{C}_o$ has NW and SE boundaries with asymptotic slopes $\rho(p)<-1$ and $1/\rho(p)>-1$ respectively.  In particular for each such $p$, the assumptions of Theorem \ref{thm:direction_trans} hold with $v=(1,1)$.  
\Qed

We believe that a similar argument shows directional transience in the direction $\uparrow$ for any model $\mc{G}_o\sim(\SWE, \uparrow)$, provided 
$\mu(\mc{G}_o=\uparrow)>p_c^{\smallFSOSP}$, where the latter critical percolation threshold is defined in 
\cite{HS_DRE1}. But we have not checked this in detail.

\section{The regeneration structure}
\label{sec:orthog01}
Assuming ellipticity, Theorems \ref{thm:dtransv} and \ref{thm:zern_merk} have been proved by Sznitman and Zerner \cite{SZ99}, Zerner \cite{Zern02} and Zerner and Merkl \cite{MZ01,Zern07}, however some parts of their arguments don't actually rely on ellipticity.  In this section we prove Theorems \ref{thm:dtransv} and \ref{thm:zern_merk} (i.e.~without any ellipticity assumption), identifying the key elements necessary to make the regeneration arguments work in general.  We believe strongly that these extensions deserve a proper and careful treatment.  
In particular, the proofs of Zerner and coauthors rely on Kalikow's 0-1 law, and our main task is to reprove this 0-1 law with no ellipticity assumption, i.e.~to prove Theorem \ref{thm:A+A-}. Since we will follow the basic outline of the arguments referenced above, it is worth isolating the novel features of what follows. There are basically two:
\begin{itemize}
\item Replacing ``ellipticity at all times'' by ``ellipticity at record times''. In other words, basing the argument on having positive probability of a fixed sequence of steps following time $t$, when $t$ is only a record time rather than an arbitrary time. 
\item Identifying (\ref{standinghypothesis}) and (\ref{eq:all_ei}) as the appropriate conditions to replace uniform ellipticity, in the sense that they give us ``ellipticity at record times''.
\end{itemize}

For fixed $d\ge 2$ we define a slab to be a region between any two parallel $d-1$ dimensional hyperplanes in $\R^{d}$.  Let $\mc{H}_{e_j}$ be the set of slabs $\mc{S}$ for which there exists a constant $H=H(\mc{S})\in \N$ such that $(\mc{S}+He_j)\cap \mc{S}=\varnothing$.  This is the set of slabs with finite width in direction $e_j$.  Note that for every $d$-dimensional slab $\mc{S}$, there exists some $j\in \{1,2,\dots, d\}$ such that $\mc{S}\in \mc{H}_{e_j}$. Moreover, if $v$ is a normal to the hyperplanes defining the slab, then $\mc{S}$ has finite width in the direction $u$ $\Leftrightarrow$ $u\cdot v\neq 0$.

In preparation for proving Theorem \ref{thm:A+A-}, we introduce the following condition, which (up to reflections of the directions) says that the RWRE is truly $d$-dimensional:
\begin{equation}
\label{eq:all_ei}
\boxed{\mu(e_i \in \mc{G}_o)>0\quad  \text{for }i=1,\dots,d.}
\end{equation}
Roughly speaking, we use the conditions (\ref{eq:all_ei}) and (\ref{standinghypothesis}) to replace the assumption of ellipticity in the regeneration arguments (that date back to \cite{K81} in the uniformly elliptic setting), and then prove the theorem both when these conditions hold, and when they fail.

\begin{LEM}
\label{cor:directional transience}
Assume (\ref{standinghypothesis}) and (\ref{eq:all_ei}).  Then for each unit vector $v$, $P$-almost surely,
\[\liminf_{n\ra \infty}X_n\cdot v\in \{-\infty, +\infty\}.\]
\end{LEM}
\proof Let $\mc{A}=\{e\in \mc{E}: \mu(e\in \mc{G}_o)>0\}$, $\mc{B}_+=\{e\in \mc{E}:e\cdot v>0\}$,
and $\mc{B}_-=\{e\in \mc{E}:e\cdot v<0\}=-\mc{B}_+$.  Note that $\{e_1,\dots, e_d\}\subset \mc{A}$ by (\ref{eq:all_ei}).  If $\mc{B}_-\subset \mc{A}^c$  then the random walk is transient in every direction in $\mc{B}_+\cap \mc{A}$
(e.g.~by Theorem \ref{thm:direction_trans}), hence it is also transient in direction $v$. 
So assume there exists some $e_-\in \mc{B}_-\cap \mc{A}$.
Without loss of generality we can assume that $e_-\in \{\pm e_1\}$. 

Assume that the claim of the lemma is false, i.e.~there exists $r\in \R$ such that $\liminf_{n \ra \infty} X_n \cdot v=r$.  Then the slab $\mc{S}=\{x\in \re^d:r-1\le x\cdot v \le r+1\}$ is visited infinitely often.   
Suppose that there are actually infinitely many sites in $\mc{S}$ that are visited. Then the set of sites of
$\mc{S}$ visited is unbounded in the direction of at least one vector $e \in U_1=\{\pm e_j:j\ne 1\}$.  Let $T_1=1$, and let
$$
T_{k+1}=\inf\{n>T_k:X_n\in\mc{S}\text{ and }\exists e\in U_1\text{ s.t. } X_n\cdot e>
\max\{X_m\cdot e:m<n, X_m\in\mc{S}\}\}
$$
be the times the walk reaches a new record level within $\mc{S}$. These are all finite, and by definition, the sites $\mc{S}\cap (X_{T_k}+\Z e_-)$
were not explored prior to time $T_k$. Let $H<\infty$ be the width of $\mc{S}$ in the direction $e_1$, defined above. With
probability at least $\epsilon^{H+1}$ there is an admissible path, just using $e_-$ arrows, of length at most $H+1$, that connects
$X_{T_k}$ to a point $z$ outside $\mc{S}$, with $z\cdot v\le r-1$. So with probability at least $\epsilon^{2(H+1)}$, the random walk follows this
path and exits $\mc{S}$ in at most $H+1$ steps. This must therefore occur almost surely, for infinitely many $k$. It follows that $\liminf_{n \ra \infty} X_n \cdot v\le r-1$, which is a contradiction.

Thus in fact there are only finitely many points of $\mc{S}$ that get visited. Therefore at least one point gets visited infinitely often. Let $\mc{R}$ be the set of sites in $\mc{S}$ that are visited infinitely often. We conclude that  $\mc{R}$ is finite but non-empty.  We may therefore choose an element $x$ of $\mc{R}$ such that $x\cdot v$ is minimal.  Since every $y \in \mc{C}_x$ is also visited infinitely often (by the argument for (i) of Lemma \ref{lem:stuck}), we must have that $(y-x)\cdot v\ge 0$ for each $y\in \mc{C}_x$, i.e.~$x\in {\bf C}_{v}$.  Since this happens with positive probability, we in fact have $\nu(o\in {\bf C}_v)>0$.  The proof of Theorem \ref{thm:direction_trans} shows that if $\mu(e_1\in\mc{G}_o)>0$ and (\ref{standinghypothesis}) holds, then on the event that $o\in{\bf C}_v$ we have $\liminf X_n\cdot v=\infty$ (a.s.), which contradicts the definition of $r$.   \Qed

\medskip
Note that we can similarly get that $\liminf X_n \cdot e_1\in \{-\infty,+\infty\}$, $P$-almost surely under the weaker assumptions (\ref{standinghypothesis}) and $\mu(e_1\in \mc{G}_o)>0$.

For any RWRE, and fixed $\ell\in \re^d \setminus o$, recall that $A_{+}=A_+(\ell)$ is the event that the walk is transient in direction $\ell$, i.e.~$A_{+}=\{X_n\cdot \ell\ra+\infty\}$, and $A_-=A_-(\ell)=\{X_n\cdot \ell\ra-\infty\}$. 
Let $O=(A_{+}\cup A_{-})^c$ and let $O_m$ be the event that both $X_n\cdot \ell\le m$ and $X_n\cdot \ell\ge m$ infinitely often.  
For $k\ge 0$ let $T_k=\inf\{n:X_n\cdot \ell\ge k\}$.
\begin{LEM}
\label{lem:Oinv}
Assume (\ref{standinghypothesis}) and (\ref{eq:all_ei}). Then $P(O)=P(\cap_{n \in \Z}O_n)$.
\end{LEM}
\proof Firstly note that $\cap_{n \in \Z}O_n\subset \cup_{n \in \Z}O_n\subset O$ so in particular if $P(O)=0$ then the result is trivial.  By Lemma \ref{cor:directional transience} applied to both $\ell$ and $-\ell$, $P$-almost surely on the event $O$ we have $\liminf X_n\cdot \ell=-\infty$ and $\limsup X_n\cdot \ell=+\infty$.  This implies that $\cap_{n \in \Z}O_n$ occurs almost surely on the event $O$ as required.
\Qed



\begin{LEM}
\label{lem:O2}
If $P(A_+)>0$ then 
\begin{equation}
P(A_+\cap \{X_n\cdot \ell\ge 0\,\forall n\ge 0\})>0.
\label{eqn:stayright}
\end{equation}
If also (\ref{standinghypothesis}) and (\ref{eq:all_ei}) then  $P(O)=0$.
\end{LEM}
\proof  Suppose that $P(A_+)>0$.  Let $\mP_{\omega,y}$ denote the quenched
law of the RW in the environment $\omega$, starting from $X_0=y$, and recall that $P_y$ is the corresponding annealed law. For any $x$ and $n$, let $A(x,n)$ be the event $A_+\cap \{X_n=x \text{ and } X_k\cdot \ell \ge x\cdot \ell\,\forall k\ge n\}$. Since $P(A_+)>0$, we must have $P(A(x,n))>0$ for some $x$ and $n$. By translation invariance, $P_{-x}(A(o,n))>0$. By the Markov property, $\mP_{\omega,-x}(A(o,n)\mid X_n=o)=\mP_{\omega,o}(A(o,0))$ for every environment $\omega$. Therefore $\mP_{\omega,-x}(A(o,n))\le\mP_{\omega,o}(A(o,0))$. Integrating out $\omega$ we get $0<P_{-x}(A(o,n))\le P(A(o,0))$, which proves (\ref{eqn:stayright}).  

To prove that $P(O)=0$, note that by Lemma \ref{lem:Oinv}, almost surely on the event $O$, all $O_m$ occur and therefore $\limsup X_n\cdot \ell=\infty$.  It is therefore sufficient to show that almost surely $O$ does not occur on the event $\limsup X_n\cdot \ell=\infty$, under the assumptions of the lemma.  
Let $\delta=P(X_k\cdot \ell\ge 0\,\forall k\ge 0)\ge P(A(o,0))>0$.
Let $T_0=0$. Given $T_k<\infty$, let $D_k=\inf\{n>T_k:X_n\cdot \ell < X_{T_k}\cdot \ell\}$. If $D_k<\infty$, let $M_k=\sup\{X_n\cdot \ell: n\le D_k\}$, and let $T_{k+1}=\inf\{n>D_k: X_n\cdot \ell\ge M_k+1\}$. Let $K$ be the first value of $k$ such that $T_k=\infty$ or $D_k=\infty$. 

At  time $T_k<\infty$ the walker has not explored any of the environment in direction $\ell$ from $X_{T_k}$, so $P(D_k=\infty\mid T_k<\infty)=\delta>0$. Therefore on the event that $\limsup X_n\cdot \ell=\infty$, we have repeated independent trials, and so eventually will have a $k$ with $D_k=\infty$, i.e.~$K<\infty$ a.s.~on $\{\limsup X_n\cdot \ell=\infty\}$. But when $K<\infty$, some $O_m$ does not occur, so by Lemma \ref{lem:Oinv} neither does $O$. 
\Qed

\noindent {\em Proof of Theorem \ref{thm:A+A-}.}
If (\ref{standinghypothesis}) fails, then by Lemma \ref{lem:stuck}, the random walk a.s.~only visits finitely many sites. Thus $P(O)=1$. 
So assume (\ref{standinghypothesis}).  Let $\mc{B}=\{i=1,\dots,d:\mu(e_i\in \mc{G}_o)+\mu(-e_i\in \mc{G}_o)>0\}$, and write
\[\ell =\sum_{i \in \mc{B}} \ell^{[i]} e_i+\sum_{i \notin \mc{B}} \ell^{[i]} e_i=\ell_{\mc{B}}+\sum_{i \notin \mc{B}} \ell^{[i]} e_i.\]
If $\ell_{\mc{B}}=0$ then $X_n\cdot \ell=0$ for all $n$ almost surely so that $P(O)=1$.  Otherwise $\ell_{\mc{B}}\ne 0$ and since $X_n \cdot \sum_{i \notin \mc{B}} \ell^{[i]} e_i=0$ for all $n$, we have that $O(\ell)\iff O(\ell_{\mc{B}})$, $A_+(\ell)\iff A_+(\ell_{\mc{B}})$ and $A_-(\ell)\iff A_-(\ell_{\mc{B}})$.  This has reduced the problem to a $|\mc{B}|$-dimensional one, so without loss of generality we may assume that $|\mc{B}|=d$.  Then by considering reflections of the axes, we may further assume that (\ref{eq:all_ei}) holds.  In this case, if $P(A_+)>0$ then $P(O)=0$ by Lemma \ref{lem:O2}. Similarly $P(A_-)>0\Rightarrow P(O)=0$ by applying Lemma \ref{lem:O2} to $-\ell$.   The result then follows from the fact that $P(O)+P(A_+\cup A_-)=1$.
\Qed

\noindent {\em Proof of Theorem \ref{thm:dtransv}.}
We have the same hierarchy of possibilities as in Theorem \ref{thm:A+A-}. If (\ref{standinghypothesis}) fails, then the random walk a.s. only visits finitely many sites, and $n^{-1} X_n\to 0$ a.s.
Likewise, if all symmetric versions of (\ref{eq:all_ei}) fail then the problem reduces to a lower dimensional one (where some symmetric version of (\ref{eq:all_ei}) with smaller $d$ does hold).  Therefore without loss of generality we assume that both (\ref{standinghypothesis}) and (\ref{eq:all_ei}) hold.
By Theorem \ref{thm:A+A-} there are two cases to consider, namely $P(A_{+} \cup A_{-})=1$ and $P(A_{+} \cup A_{-})=0$.

The former is addressed in Theorem 3.2.2 of \cite{Zeit04}, for i.i.d.~uniformly elliptic environments, drawing on ideas that go back to \cite{K81}, together with contributions of Zerner and \cite{SZ99}. As pointed out in \cite{Zern02}, the proof does not actually require uniform ellipticity, but works simply assuming ellipticity. In fact, even ellipticity is used only to obtain (\ref{eqn:stayright}). In other words, the argument of \cite{Zeit04} applies to IID environments satisfying (\ref{eqn:stayright}) and $P_o(A_{+} \cup A_{-})=1$. In particular, this proves the theorem when $P(A_{+} \cup A_{-})=1$.

For completeness, we sketch the argument. Adopting notation from the proof of Lemma \ref{lem:O2}, let $\tau_1=D_K$. On the event $A_{+}$, $\tau_1$ acts as a regeneration time, so conditional on $A_{+}$, the process $\hat X_n=X_{\tau_1+n}-X_{\tau_1}$ and the environment $\hat \omega_x=\omega_{x+X_{\tau_1}}$ (for $x\cdot \ell\ge 0$) are independent of the environment and walk observed up to time $\tau_1$. This allows one to construct additional regeneration times $\tau_{1}<\tau_2<\dots$ such that (conditional on $A_{+})$ the $X_{(\tau_k+n)\land \tau_{k+1}}-X_{\tau_k}$ are IID segments of path. If $E[\tau_1]<\infty$, the strong law of large numbers now implies the existence of a deterministic speed $v_{+}$ on the event $A_{+}$. If $E[\tau_1]=\infty$, one appeals to a calculation \cite[Lemma 3.2.5]{Zeit04} due to Zerner,
which shows that 
\begin{equation}
E[(X_{\tau_{k+1}}-X_{\tau_{k}})\cdot \ell|A_{+}]\le \frac{C}{P(D_0=\infty)}<\infty. 
\label{eqn:renewaltheory}
\end{equation}
This estimate is enough to give, by the law of large numbers again, that the speed $v_{+}$ exists on $A_{+}$ and $=0$.  Note that the proof of \cite[Lemma 3.2.5]{Zeit04} is presented with $\ell=e_1$, in which case the left side of \eqref{eqn:renewaltheory} is actually shown to equal $1/P(D_0=\infty)$. With general $\ell$ we have only been able to verify the inequality given in \eqref{eqn:renewaltheory}, but that certainly suffices for our purposes. 

It remains to show the case $P(A_{+} \cup A_{-})=0$. This is addressed in Theorem 1 of \cite{Zern02}, again assuming ellipticity. But the proof carries over verbatim if ellipticity is replaced by the following weaker condition:
\begin{multline*}
\text{If $\{x\in\Z^d: a\le x\cdot \ell\le b\}$ is visited by $X_n$ infinitely often, }\\
\text{ then there a.s. exist $n, m$ with $X_n\cdot \ell<a$ and $X_m\cdot \ell>b$.}
\end{multline*}
The latter property holds in our setting, by applying Lemma \ref{cor:directional transience} to $\ell$ and $-\ell$. 
 \Qed
 

\begin{COR}
\label{cor:dtransv}
Assume that $P(A_\ell)\in \{0,1\}$ for each $\ell \in \{ e_1,\dots, e_d\}$. Then there exists a deterministic $v\in \re^d$ such that
\[\lim_{n\ra \infty}\frac{X_n}{n}=v.\]
\end{COR}
\proof
As in Corollary 2 of \cite{Zern02}, apply Theorem \ref{thm:dtransv} to each coordinate direction. Note that by Theorem \ref{thm:A+A-}, $P(A_{e_k})\in \{0,1\}\Rightarrow P(A_{-e_k})\in \{0,1\}$.
\Qed

\noindent {\em Proof of Theorem \ref{thm:zern_merk}.}  Zerner and Merkl prove this (Theorem 1 of \cite{MZ01}, see also \cite{Zern07}) assuming ellipticity. In fact, the proof is valid under the following conditions: $P(A_{+}\cup A_{-})\in\{0,1\}$; $P(A_{+})>0\Rightarrow P(X_n\cdot \ell\ge 0\,\forall n)>0$; 
$P(A_{-})>0\Rightarrow P(X_n\cdot \ell\le 0\,\forall n)>0$. The first property holds in our setting, by Theorem \ref{thm:A+A-}.  The second and third properties hold by Lemma \ref{lem:O2}.
\Qed


\begin{COR}
\label{cor:speedsin2d}
When $d=2$ there exists a deterministic $v\in \re^2$ such that
\[\lim_{n\ra \infty}\frac{X_n}{n}=v.\]
\end{COR}
\proof
Corollary \ref{cor:dtransv} and Theorem \ref{thm:zern_merk}.
\Qed

\section{Properties obtained by coupling}
\label{sec:coupling}
In this section we use coupling methods to prove a number of results, beginning with the monotonicity result of Theorem \ref{thm:main1}.

\subsection{Monotonicity}
\label{sub:coupling_mono}
For a RWRE $X$ in an environment taking at most countably many values $\gamma^{\sss i}$, $i\in \N$, we define
$E^{\sss i}=\{e\in \mc{E}: \gamma^{\sss i}(e)>0\}$ and 
\[N_n^{\sss i}=\{0\le m<n:\omega_{X_m}=\gamma^{\sss i}\}.\]
We also let $u_i=\sum_{e\in\mc{E}}\gamma^{\sss i}(e)e$ be the local drift of environment $\gamma^{\sss i}$.

For a $2$-valued environment we also write $N_n=N_n^{\sss 1}$.  Note that in this case, for almost every environment, $N_n^{\sss i}\ra \infty$ almost surely for $i=1,2$ even if the walker gets stuck on $\ge 2$ sites.

\begin{THM}
\label{thm:other_mono}
For any $2$-valued model $(\gamma^{\sss 1},\gamma^{\sss 2})$ with $\gamma^{\sss 1},\gamma^{\sss 2}\ne \emptyset$, there exists a coupling under which for all $0\le p<p'\le 1$ the following hold:
\begin{itemize}
\item[(i)]  $N_n[p']\ge N_n[p]$ almost surely, 
\item[(ii)] for every $e$ such that $e\in E^{\sss 1}\cap (E^{\sss 2})^c$ and $-e\notin E^{\sss 1}$, $X_n[p']\cdot e\ge X_n[p]\cdot e$ a.s.~for every $n\ge 0$, 
\item[(iii)] $B_n=\{X_n^{[1]}[p']\ge X_n^{[1]}[p]\}$ occurs for infinitely many $n$, almost surely if $u_1^{[1]}>u_2^{[1]}$.  
\item[(iv)] Let $V[p]=\lim_{n\to\infty} n^{-1}X_n[p]$ (which exists a.s. by Theorem \ref{thm:dtransv}, but could in principle be random). Then for any $u\ne o$,  
\begin{itemize}
\item[(a)] $V[p']\cdot u\ge  V[p]\cdot u$ a.s. if  $(u_1-u_2)\cdot u\ge 0$, and
\item[(b)] $V[p']\cdot u\le  V[p]\cdot u$ a.s. if $(u_1-u_2)\cdot u\le 0$.
\end{itemize}

\end{itemize}
\end{THM}
\proof
Let $\{U_x\}_{x\in \Z^d}$, $\{Y^{\sss 1}_n\}_{n\in \N}$, and $\{Y^{\sss 2}_n\}_{n\in \N}$ be independent random variables with distributions $U[0,1]$, $\gamma^{\sss 1}$, and $\gamma^{\sss 2}$ under $\mP$ respectively.  Define $\omega[p]=(\omega_x[p])_{x\in \Z^d}$ by 
\[\omega_x[p]=\begin{cases}
\gamma^{\sss 1} & \text{ if }U_x<p\\
\gamma^{\sss 2} & \text{ otherwise. }
\end{cases}\]  

Set $X_0[p]=0$ and given $X_0[p],\dots,X_n[p]$ define,
\[X_{n+1}[p]-X_n[p]=\begin{cases}
Y^{\sss 1}_k, &\text{ if }\omega_{X_n}[p]=\gamma^{\sss 1}, \text{ and } N_n[p]=k-1\\
Y^{\sss 2}_k, &\text{ if }\omega_{X_n}[p]=\gamma^{\sss 2}, \text{ and } n-N_n[p]=k-1.
\end{cases}\]
One can easily check that $X[p]$ is a RWRE in environment $\omega[p]$.  

For the first claim, let $p'>p$.  Define 
$T_1=\inf\{n\ge 1:N_n[p]\ne N_n[p']\}=\inf\{n\ge 1:N_n[p]<N_n[p']\}$ and $T^*_1=\inf\{n\ge T_1+1:N_n[p']=N_n[p]\}$.  
The claim certainly holds up to time $T^*_1$.  Therefore if $T_1=\infty$, or if $T_1<\infty$ and $T^*_1=\infty$, then there is nothing to prove.  So assume $T_1, T^*_1<\infty$.  We have $N_{T^*_1}[p']=N_{T^*_1}[p]$.  Under the given coupling $X_{T^*_1}[p']$ and $X_{T^*_1}[p]$ have therefore taken exactly the same number of steps in each direction so that $X_{T^*_1}[p']=X_{T^*_1}[p]$, and the walks are recoupled.  We can now repeat the above argument with 
$T_2=\inf\{n>T^*_1\ge 1:N_n[p]\ne N_n[p']\}=\inf\{n\ge 1:N_n[p]<N_n[p']\}$ and $T^*_2=\inf\{n\ge T_2+1:N_n[p']=N_n[p]\}$, etc. to get (i). 

Suppose $e \in E^{\sss 1}\cap (E^{\sss 2})^c$ but $-e\notin E^{\sss 1}$. Then the number of $e$-steps taken by the walk $X[p]$ up to time $n$ is $\#\{k\le N_n[p]: Y_k^{\sss 1}=e\}$.  The number of $-e$ steps taken is $\#\{k\le n-N_n[p]: Y_k^{\sss 2}=-e\}$.  The second claim now follows since $N_n[p]$ is increasing in $p$ for each $n$ under this coupling.

To prove (iii), let $W_j^{\sss r}=\mathbbm{1}_{\{Y^{\sss r}_j=e_1\}}-\mathbbm{1}_{\{Y^{\sss r}_j=-e_1\}}$. Observe that $E[W_j^{\sss 1}]=u_1^{[1]}>u_2^{[1]}=E[W_j^{\sss 2}]$ and 
\begin{align*}
X_n^{[1]}[p']-X_n^{[1]}[p]=\sum_{j=N_n[p]+1}^{N_n[p']}W_j^{\sss 1} - \sum_{j=n-N_n[p']+1}^{n-N_n[p]}W_j^{\sss 2},
\end{align*}
where each of these sums contains $N_n[p']-N_n[p]\ge 0$ elements.  Define $S_0=0$ and for $i\in \N$,
\begin{align*}
R_i=&\inf\{n>S_{i-1}:N_n[p']>N_n[p]\} \quad \text{and} \quad S_i=S_i^*\vee S_i^{**}, \quad \text{where}\\
S_i^*=&\inf\{n>R_i:N_n[p]= N_{R_i}[p']\} \quad \text{and} \quad S_i^{**}=\inf\{n>R_i:n-N_n[p']= R_i-N_{R_i}[p]\}.
\end{align*}
Note that since $N_n[p']\ge N_n[p]$ for all $n$ and the walks are at the same location whenever $N_n[p']=N_n[p]$ we also have that $R_i<\infty$ for each $i$ almost surely (while the walks stay together, consider the times at which they reach a new record level).
Then also $S_i<\infty$ since both $N_n$ and $n-N_n$ are increasing to $+\infty$ 
(this property does not depend on $p\in (0,1)$).   

By definition, $N_{S_i}[p']\ge N_{S_i}[p]\ge N_{R_i}[p']\ge N_{S_{i-1}}[p']$ so the intervals  $\big\{\big[N_{S_i}[p]+1,N_{S_i}[p']\big]\big\}_{i\ge 1}$ are disjoint (possibly empty).  Similarly $S_i-N_{S_i}[p]\ge S_i-N_{S_i}[p']\ge R_i-N_{R_i}[p]\ge S_{i-1}-N_{S_{i-1}}[p]$ so the intervals  $\big\{\big[S_i-N_{S_i}[p']+1,S_i-N_{S_i}[p]\big]\big\}_{i\ge 1}$ are disjoint.
It follows that the random variables
\[\Delta_i\equiv X_{S_i}^{[1]}[p']-X_{S_i}^{[1]}[p]=\sum_{j=N_{S_i}[p]+1}^{N_{S_i}[p']}W_j^{\sss 1} - \sum_{j={S_i}-N_{S_i}[p']+1}^{{S_i}-N_{S_i}[p]}W_j^{\sss 2}, \quad i\in \N\]
are sums of random variables that are independent of the random variables summed in $\{\Delta_{j}\}_{j<i}$.  Note that the lengths of the intervals we sum over need not be independent.

Fix $k\in \N$ and let $A_i=\{\Delta_i\ge 0\}$.  Then either
\begin{align}
\label{e:condition}
P\left(\ccap{i=k}{\infty}A_i^c\right)=P(A_k^c)\prod_{j=k+1}^{\infty}P\left(A_j^c \Big|\ccap{i=k}{j-1}A_i^c\right),
\end{align}
or the left hand side of \eqref{e:condition} is zero.  By the law of large numbers, and using the fact that $E[W_j^{\sss 1}]>E[W_j^{\sss 2}]$, there is some $\varepsilon>0$ such that $P\big(\ccap{m=1}{\infty} \{\sum_{i=1}^mW_j^{\sss 1}\ge \sum_{i=1}^mW_j^{\sss 2}\}\big)>\varepsilon$.  Thus 
\[P\left(A_j^c \Big|\ccap{i=k}{j-1}A_i^c\right)<1-\vep\]
for every $j$. It follows then that \eqref{e:condition} is equal to zero and hence 
\[P(B_n \text{ i.o.})=\lim_{k\ra \infty}P(\ccup{i=k}{\infty}B_i)\ge  \lim_{k\ra \infty} P(\ccup{i=k}{\infty}A_i)=1.\]

To prove (iv), note that 
\begin{align}
\label{e:XNrelation}
\frac{X_n}{n}=&\frac{N_n}{n}\frac{1}{N_n}\sum_{j=1}^{N_n}Y^{\sss 1}_j+\left(1-\frac{N_n}{n}\right)\frac{1}{n-N_n}
\sum_{j=1}^{n-N_n}Y^{\sss 2}_j
 = V^{\sss 2}_{n}+(V^{\sss 1}_{n}-V^{\sss 2}_{n})\frac{N_n}{n},
\end{align}
where $V^{\sss i}_{n}\equiv\frac{1}{N^{\sss i}_n}\sum_{j=1}^{N^{\sss i}_n}Y_j^{\sss i}\to u_i$ a.s. as $n\ra \infty$ for $i=1,2$ by the LLN and the fact that $N_n^{\sss i}\ra \infty$. 
Therefore for each $p\in [0,1]$,
\begin{equation}
V[p]\cdot u=u_2\cdot u+\lim_{n\ra \infty}(u_1-u_2)\cdot u\frac{N_n[p]}{n}.
\label{velocityrelation}
\end{equation}
The limits on the right therefore exist. Moreover, on this probability space $N_n[p]$ is almost surely increasing in $p$, so the quantities $V[p]\cdot u$ are a.s. increasing in $p$ if $(u_1-u_2)\cdot u\ge 0$, and a.s. decreasing in $p$ if $(u_1-u_2)\cdot u\le 0$.  
\Qed

Note that as stated, (iv) of Theorem \ref{thm:other_mono} is a statement about coupled random variables. But when $d=2$ we know that deterministic speeds exist, so this becomes simply a statement about monotonicity of those speeds.\medskip

Theorem \ref{thm:other_mono} applied to each of the 2-valued RWRE models $(\NWE, \downarrow)$, $(\NWE, \WE)$ and $(\NWE, \SWE)$ shows that there exists a coupling under which $X_n(p)\cdot(0,1)$ is almost surely increasing in $p$ for all $n$.  Applied to the model $(\NE, \SW)$ the theorem gives a coupling under which $X_n[p]\cdot (1,1)$ is almost surely increasing in $p$ for all $n$.  This in turn implies that for $p\ge \hlf$ the probability that the model is transient in the direction $-(1,1)$ is at most $\hlf$, and by Theorem \ref{thm:zern_merk}, it must be 0.  

Theorem \ref{thm:other_mono}(ii) implies that in the uniform RWRE on ($\NSE$, $\NSW$), $X_n^{[1]}[p]$ is monotone in $p$, but it is easily checked that for the given coupling, this statement may fail for the uniform RWRE on ($\NSE$, $\NSEWalt\,\,\,$). On the other hand, Theorem \ref{thm:other_mono}(iii) establishes monotonicity for the latter model in a weaker sense.

In 2 dimensions, we know that a deterministic limit $v[p]=\lim_{n\to\infty}n^{-1}X_n$ always exists. This  
implies the following simple linear relationship between $v^{[1]}[p]$ and $v^{[2]}[p]$ that is independent of $p$ (see e.g.~Table \ref{tab:walks}).
\begin{COR}
\label{cor:vrelation}
For any $2$-valued model with $d=2$,
$$
(v^{[2]}[p]-u_2^{[2]})(u_1^{[1]}-u_2^{[1]})=(v^{[1]}[p]-u_2^{[1]})(u_1^{[2]}-u_2^{[2]}).
$$
\end{COR}
\proof  By \eqref{velocityrelation}, both sides equal $\lim_{n\to\infty}(u_1^{[1]}-u_2^{[1]})(u_1^{[2]}-u_2^{[2]})\frac{N_n[p]}{n}$.\Qed

\begin{COR}
\label{cor:asymptoticfrequency}
For any $2$-valued model, $\lim_{n\to\infty}n^{-1}N_n$ exists a.s. (but may be random).
\end{COR}
\proof 
If $u_1\neq u_2$ then the claim holds by (\ref{velocityrelation}). 
If $u_1=u_2\neq o$, then $n^{-1}X_n\to u_1$. Thus regeneration times will exist, as in Theorem 
\ref{thm:dtransv}, which implies the convergence of $n^{-1}N_n$ by the law of large numbers. Finally, if $u_1=u_2=o$ then the random environment is {\it balanced}, in the sense of Berger and Deuschel \cite{BD11}. In that paper they establish the existence of an ergodic measure $Q$, absolutely continuous with respect to $P$, and invariant for the environment as viewed from the particle. This again is sufficient to imply convergence of $n^{-1}N_n$.
\Qed

The above result can be generalized beyond the 2-valued case.  When (\ref{standinghypothesis}) does not hold (see also Theorem \ref{thm:stuck}) and the environment is at least 4-valued, it is possible that there are two or more distinct arrangements of sites on which the walk can get stuck.  In this case the asymptotic frequency $\lim_{n\to\infty}n^{-1}N^{\sss i}_n$ exists a.s., but will depend on which trapping configuration the walk finds.
2-valued models can get trapped for long periods in a finite set of sites (forever if $\mc{G}_o\in \{e,-e\}$ a.s.~for some $e\in \mc{E}$), but the configurations which accomplish this are essentially trivial. 

\medskip
Theorem \ref{thm:other_mono} gives monotonicity results that apply to all 2-valued i.i.d.~random environments.  
The logical next step is to investigate these questions when the local environment can take more than 2 values.  However, as soon as we have more than 2 possible environments, we also have multiple different notions of monotonicity.  Some possibilities are given below. In each case the environments are $\gamma^1, \dots, \gamma^k$, with probabilities $p_1, \dots, p_k$, and we consider velocities $v^{[i]}=\lim_{n\to\infty}n^{-1}X^{[i]}_n$ and frequencies $\pi^{\sss i}=\lim_{n\to\infty}n^{-1}N^{\sss i}_n$. 
\begin{itemize}
\item[(i)] Monotonicity of $v^{[1]}$ as we $\uparrow p_1$ and $\downarrow p_2$, holding fixed $p_3, \dots, p_k$. 
\item[(ii)] Monotonicity of $\pi^{\sss 1}$ as we $\uparrow p_1$ and $\downarrow p_2$, holding fixed $p_3, \dots, p_k$. 
\item[(iii)] Monotonicity of the ratio $\pi^{\sss 1}/\pi^{\sss 2}$ as we $\uparrow p_1$ and $\downarrow p_2$, holding fixed $p_3, \dots, p_k$. 
\item[(iv)] Monotonicity of $v^{[1]}$ or $\pi^{\sss 1}$ as we vary $p_1$, holding fixed the {\it relative} sizes of $p_2,\dots,p_k$. 
\item[(v)] Monotonicity of $v^{[1]}$ or $\pi^{\sss 1}$ as the weight to $e_1$ within $\gamma^{\sss 1}$ increases, holding all probabilities $p_1,\dots,p_k$ fixed.  
\end{itemize}
We will give counterexamples to some of these below. In each case, $(u_1^{[1]}-u_2^{[1]})\cdot e_1>0$, so monotonicity would hold in the 2-valued case. See \cite{HSun10} for a monotonicity result under additional strong assumptions that allows for arbitrarily many values of the environment. See \cite{MadrasT} for an early example of non-monotonicity in random walk models.

\begin{figure}
\includegraphics[scale=.4]{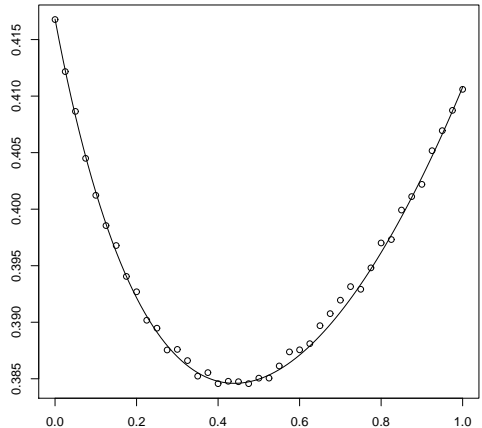}
\caption{The speed $v^{[1]}(q)$  for $q\in(0,1)$ in Example \ref{exa:3valuedspeed}, with $\alpha=.95$ and $p_3
=.03$ fixed,  
and $p_1=.97q$ and $p_2=.97(1-q)$. The curve is the exact 
value, while each point is an estimate based on 1,000 simulations of 20,000 step random walks, carried out in R.}
\label{fig:speed_non_mono}
\end{figure}

\begin{EXA}
\label{exa:3valuedspeed} $\NEalpha\vspace{2mm} \NE\leftarrow:$
(a 3-valued counterexample to monotonicity of $v^{[1]}$ as in (i))
\end{EXA}

Let $\gamma^{\sss 1}(e_1)=\alpha=1-\gamma^{\sss 1}(e_2)$, $\gamma^{\sss 2}(e_1)=\frac12=1-\gamma^{\sss 2}(e_2)$, and $\gamma^{\sss 3}(-e_1)=1$.  Fix $p_3=p$ to be small, and vary $p_1$ and $p_2$ by setting $p_1=q(1-p)$ and $p_2=(1-q)(1-p)$, where $0\le q\le 1$. The heuristic is as follows: For $q=0$ the speed is independent of $\alpha$, and positive. For $q=1$, the speed will be greater than this value when $\alpha>\frac12$ is moderate (since $p$ is small), but declines to 0 because of trapping as $\alpha\uparrow 1$. Therefore one expects to find an alpha so that the speeds with $q=0$ and $q=1$ match, and so as long as the speed is not actually constant, it will be non-monotone at some point in between.

To actually prove non-monotonicity for this example, the most direct approach is to calculate the velocity explicitly using the regeneration structure in the direction of $e_2$ (as in Lemma \ref{lem:NE_NWspeed}).  In the notation of that result, we have that $v^{[1]}=\mE[X^{[1]}_T]/E[T]$, where $T$ is the first time $n$ that $X_n^{[2]}=1$.  
Let $\Gamma\in (0,1)$ be a random variable having the distribution of $\omega_o(e_1)$ conditional on the event $\omega_o\ne \gamma^3$.  Let $G_{\Gamma}\ge 0$ be a random variable such that conditionally on $\Gamma$, $G_{\Gamma}\sim$ Geometric$(\Gamma)$.  Let 
\begin{align*}
\eta=&\mE[\Gamma]=q\alpha+\frac{1-q}{2} \qquad \text{ and } \qquad 
\xi=\mE[G_{\Gamma}]=\frac{q\alpha}{1-\alpha}+1-q. 
\end{align*}
Then (see \cite{HS_speeds} for details)
\begin{align*}
\mE[X^{[1]}_T]=\frac{(1-p)^2\eta}{1-(1-p)\eta}-\frac{p}{1-p}   \qquad \text{ and } \qquad 
\mE[T]=p(1+2\xi)+\frac{p}{1-p}    +\frac{1-p}{1-(1-p)\eta}\left[1+2p\xi\right],
\end{align*}
One can then evaluate $v^{[1]}$ at various parameter values to find explicit examples satisfying the claim of the Lemma.  For example, the speed is non-monotone in $q$ when $p=.03$ and $\alpha=.95$ as in Figure \ref{fig:speed_non_mono}.

\begin{figure}
\includegraphics[scale=.4]{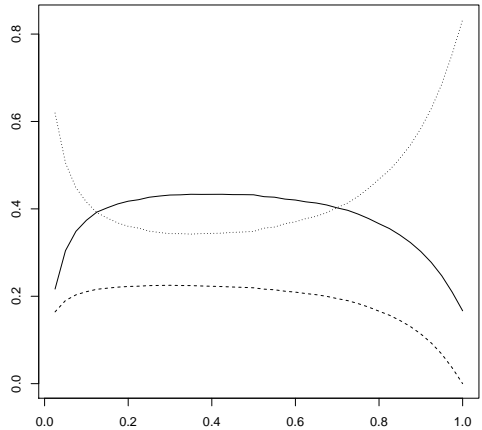}
\caption{Estimates of $\pi^{\sss i}(q)$ for $q\in (0,1)$ in Example \ref{exa:3valuedproportions}, with $p_3=.9$ and $\beta=.01$ fixed, and $p_1=.1q$ and $p_2=.1(1-q)$.
Solid, dashed and dotted lines correspond to $\pi^{\sss i}$ for $i=1,2,3$ respectively.  
Each curve is based on 3,000 repetitions of 20,000 step random walks for each of 40 values of $q$ between $0$ and $1$, carried out in R.}
\label{fig:N_nonmono}
\end{figure}

\begin{EXA}
\label{exa:3valuedproportions} $(\NE$, $\EWbeta$, $\downarrow):$
(a 3-valued counterexample to monotonicity of $\pi^{\sss 1}$ as in (ii))
\end{EXA}
Consider a 3-valued environment $\gamma^{\sss 1},\gamma^{\sss 2},\gamma^{\sss 3}$ with probabilities $p_1=q(1-p)$, $p_2=(1-q)(1-p)$, and $p_3=p$ respectively, such that
\[\gamma^{\sss 1}(e_1)=\hlf=\gamma^{\sss 1}(e_2), \quad \gamma^{\sss 2}(e_1)=\beta=1-\gamma^{\sss 2}(-e_1), \quad \gamma^{\sss 3}(-e_2)=1.\]
Here $\beta$ is small, $p_3=p$ is moderately large, and $q$ varies from 0 to 1. The heuristic is as follows: 
When $q$ is sufficiently close to $0$ there are very few traps, and the proportion of time spent at $\gamma^{\sss 1}$ sites is also close to $0$.  When $q$ is sufficiently close to $1$ there also are few traps, so since $p$ is moderately large, the proportion of time spent at $\gamma^{\sss 1}$ sites is $\ll .5$.  When $q$ is intermediate between these values, there are traps of the form 
$$
\begin{matrix}
\gamma^{\sss 3} \vspace{1mm}&\\
\gamma^{\sss 1} &\gamma^{\sss 2} 
\end{matrix}
\qquad =\qquad 
\begin{matrix}
\downarrow \vspace{1mm}&\\
\,\,\,\,\, \NE & \hspace{-2mm}\lowEWbeta
\end{matrix}
\vspace{2mm}
$$
whence the proportion of time spent at $\gamma^{\sss 1}$ sites is about $0.5$.  See Figure \ref{fig:N_nonmono} for simulations that support this heuristic, with $p=.9$ and $\beta=.01$;

\blank{
\begin{EXA}
\label{exa:4valuedspeed} $\EWbeta \NEalpha \SEalpha \rightarrow:$
(a 4-valued counterexample to monotonicity of $\pi^{\sss 1}$, with $p_2:p_3:p_4$ constant, as in (iv))
\end{EXA}
$\gamma^{\sss 1}=\EWbeta$ is as in Example \ref{exa:3valuedproportions}, and $\gamma^{\sss 2}=\NEalpha$ is as in Example \ref{exa:3valuedspeed}. Likewise $\SEalpha\vspace{2mm}$ is the $\gamma^{\sss 3}$ for which $\gamma^{\sss 3}(e_1)=\alpha=1-\gamma^{\sss 3}(-e_2)$, and $\gamma^{\sss 4}=\rightarrow$. We vary $p_1=q$, and let $p_2=p(1-q)=p_3$ and $p_4=(1-2p)(1-q)$ move proportionately, so the parameters of the model are $\alpha$, $\beta$, and $p$. 

Take $2p$ reasonably close to 1 and $\beta$ reasonably close to 0, but $\alpha$ very small. When $q=0$, $\pi^{\sss 1}$ starts out at 0. For small $q$ the most likely traps involve a single $\gamma^{\sss 1}$ so should have the form $\rightarrow \EWbeta$. Therefore $\pi^{\sss 1}$ rises quickly to be about .5; As $q$ increases further, traps involving two $\gamma^{\sss 1}$ become more likely, and therefore predominate because of the small size of $\alpha$. In other words, the walker spends most of its time in configurations of the form
$$
\begin{matrix}
\gamma^{\sss 3} \vspace{3mm}&\gamma^{\sss 1}\\
\gamma^{\sss 2} &\gamma^{\sss 1} 
\end{matrix}
\qquad =\qquad 
\begin{matrix}
\SEalpha &\highEWbeta\\
\NEalpha & \lowEWbeta
\end{matrix}
$$
Since $\alpha$ is reasonably small, this brings $\pi^{\sss 1}$ back down, before it finally rises to 1 as $q\uparrow 1$. These heuristics are supported by \dots {\red verify and insert pictures}
}
\begin{figure}
\includegraphics[scale=.4]{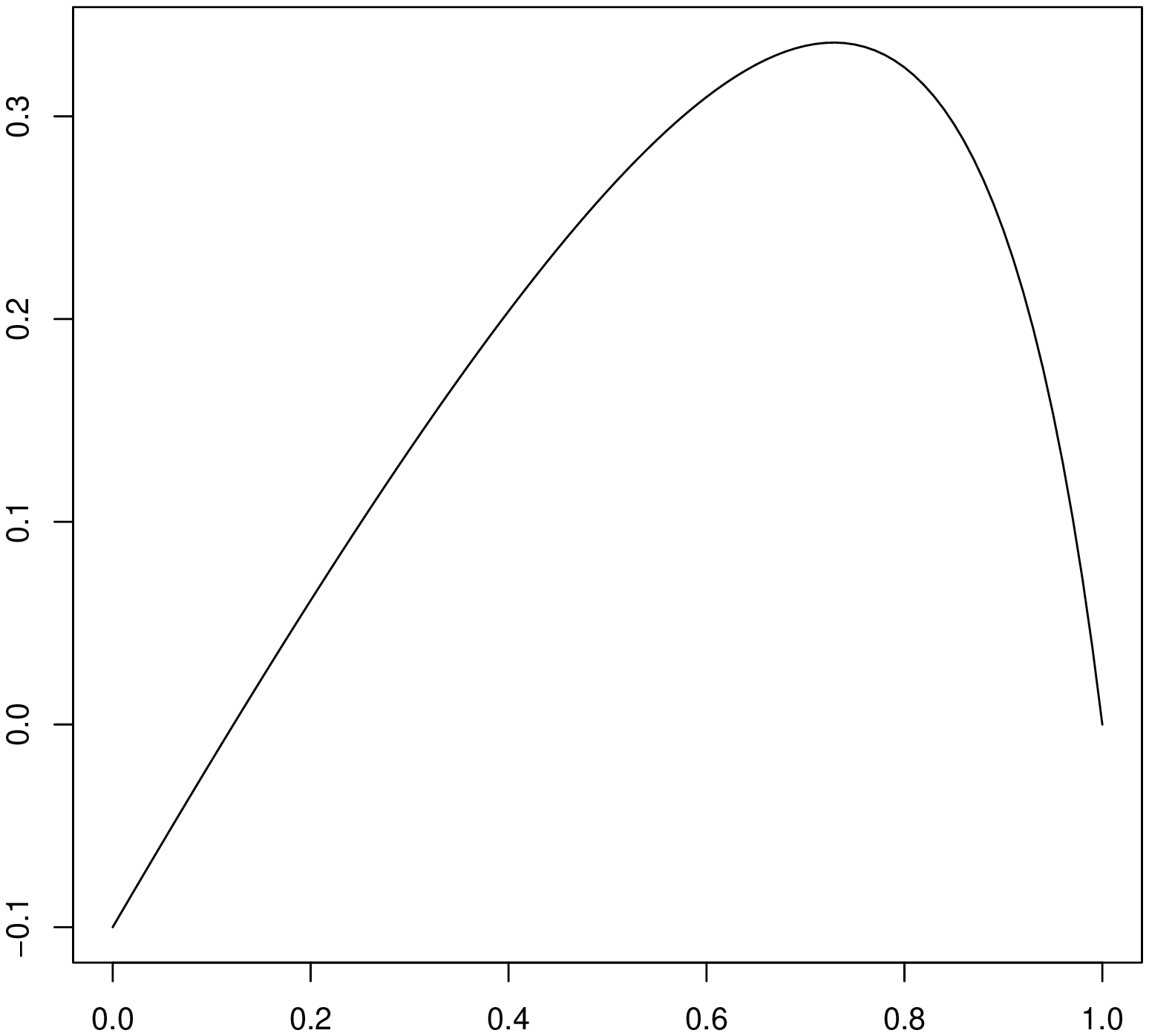}\includegraphics[scale=.4]{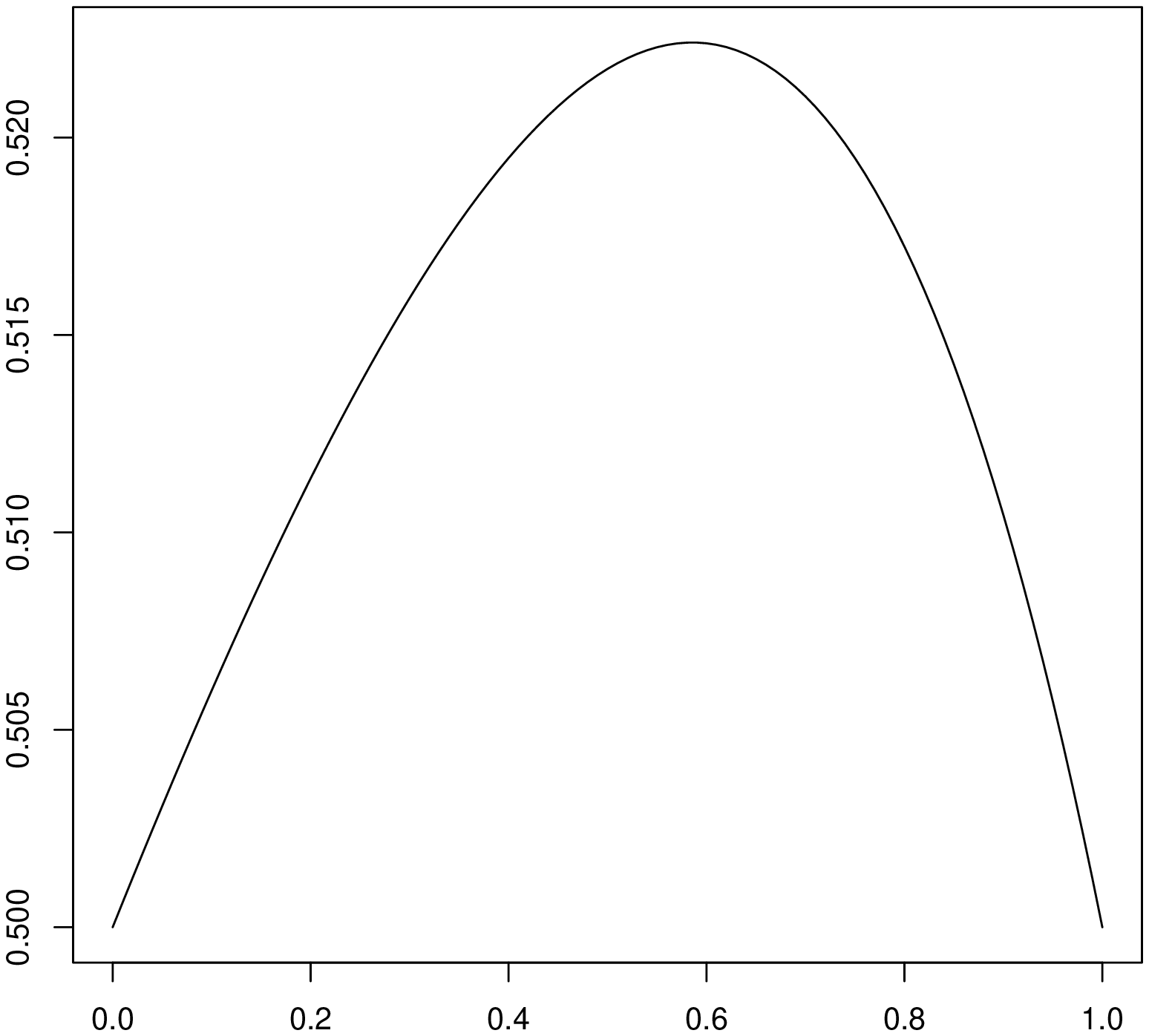}
\caption{The quantities $v^{[1]}(\alpha)$ with $p=.9$ and $\pi^{\sss 1}(\alpha)$ with $p=.5$, for $\alpha \in [0,1]$, in Example \ref{exa:2valuedspeed}.}
\label{fig:non_mono_alpha}
\end{figure}

\blank{
pi_1=function(p,alpha){1/(alpha+(1-alpha)*(p+1/p -1)+(1-p)/(1-p*alpha))}
v_1=function(p,alpha){alpha*pi_1(p,alpha)-(1-pi_1(p,alpha))}
}

\begin{EXA}
\label{exa:2valuedspeed} $\NEalpha\vspace{2mm} \leftarrow:$
(a 2-valued counterexample to monotonicity of $v^{[1]}$ in $\alpha$, as in (v))
\end{EXA}
$\gamma^{\sss 1}=\NEalpha\vspace{1mm}$ is as in Example \ref{exa:3valuedproportions}, while $\gamma^{\sss 2}$ is the environment 
$\leftarrow$, i.e.~$\gamma^{\sss 2}(-e_1)=1$. 
With $p_1=p$ and $p_2=1-p$, an elementary calculation as in 
Lemma \ref{lem:NE_NWspeed} (see \cite{HS_speeds}) shows that 
$$
\pi^{\sss 1}=\left(\alpha+(1-\alpha)\Big(p+\frac1{p}-1\Big)+\frac{1-p}{1-p\alpha}\right)^{-1}, \quad
v^{[1]}=\alpha\pi^{\sss 1}-(1-\pi^{\sss 1}).
$$
Both are non-monotone in $\alpha$ for fixed suitable choices of $p$. For example, $p=\frac12$ implies that $\pi^{\sss 1}=\frac12$ for both $\alpha= 0$ and $\alpha= 1$, without actually being constant. Likewise we have $v^{[1]}\le 0$ for both $\alpha=0$ and $\alpha=1$, but for $0<\alpha<1$, $v^{[1]}$ will be $>0$ when $p$ is large enough.  See Figure \ref{fig:non_mono_alpha}.

We have a conjectured counterexample for behaviour (iii), which we are still exploring.

\subsection{Transience}
In this subsection we begin by stating a trivial coupling criterion, which guarantees that the RWRE is transient when some related walk is transient.  We apply this criterion to prove transience results for some of our models. 
\begin{LEM}
\label{lem:coup}
Suppose that a RWRE $\{X_n\}_{n\ge 0}$ can be coupled with a random walk $X'_n$, such that for all $n,m\ge 0$, $X_n=X_m \Rightarrow X'_n=X'_m$.  Then 
\begin{itemize}
\item if $\{X'_n\}_{n\ge 0}$ is transient (almost surely) then so is $\{X_n\}_{n\ge 0}$
\item if $\{X'_n\}_{n< m}\cap \{X'_n\}_{n\ge m}=\varnothing$ then $\{X_n\}_{n< m}\cap \{X_n\}_{n\ge m}=\varnothing$ (i.e~cut-times for $X'$ are also cut-times for $X$).
\end{itemize}
\end{LEM}

A natural application of this result is the following result, that concerns the high-dimensional analogue of the uniform RWRE of Example \ref{exa:NE_SW} (which we call the {\em orthant} model). Recall that $E^i=\{e: \gamma^{\sss i}(e)>0\}$. 
\begin{COR}
Let $X=X(d,p)$ denote the uniform RWRE $(\gamma^{\sss 1},\gamma^{\sss 2})$ with $E^1=\{e_i,i=1,\dots, d\}$ and $E^2=-E^1$.  Then 
\begin{itemize}
\item when $d\ge 6$, $X$ is transient for all $p$, $P_{[p]}$-almost surely, and 
\item when $d\ge 10$, for each $p$, there exists $v[p]$ with $v^{[i]}[p]$ non-decreasing in $p$ for each $i=1,\dots,d$ such that $P_{[p]}(n^{-1}X_n \ra v[p])=1$.
\end{itemize}
\end{COR}
\proof Fix $d\ge 6$ and define a $d'$-dimensional (with $d'=\lfloor d/2\rfloor$) random walk $\{Y_n\}_{n\ge 0}$ by $Y_0=X_0=o$, and for $n\ge 1$
\begin{equation}
Y_{n}-Y_{n-1}=\begin{cases}
+e_i, \text{ if } X_n-X_{n-1}\in \{+e_{2i-1},-e_{2i}: 2i\le d\}\\
-e_i, \text{ if } X_n-X_{n-1}\in \{-e_{2i-1},+e_{2i}:2i\le d\}\\
0, \text{ otherwise}.\end{cases}
\end{equation}
Then $(Y_n^{[1]},\dots,Y_n^{[d']})= (X_n^{[1]},\dots,X_n^{[d]})A$ where $A^{[i,j]}=\mathbbm{1}_{j=2i-1}-\mathbbm{1}_{j=2i}$, i.e.
\[A^t=\begin{pmatrix}
1&-1&0&0&0&0&0&\dots &0\\
0&0&1&-1&0&0&0&\dots &0\\
0&0&0&0&1&-1&0&\dots &0\\
\vdots& & & && \ddots &\ddots & &0
\end{pmatrix}.\]
Clearly then, $Y_n=Y_m$ whenever $X_n=X_m$.  If $d$ is even, then $Y$ is a simple random walk in $d'$ dimensions.  If $d$ is odd, then $Y$ is a random walk with nearest neighbour steps but also a $\frac1d$ probability of $Y$ in place. Thus $Y$ is transient when $d'\ge 3$, so $X$ is transient when $d\ge 6$ by Lemma \ref{lem:coup}.  When $d'\ge 5$, $Y$ has well behaved cut-times. Therefore so does $X$, so it is shown in \cite{BSZ03} (see also \cite{HSun10}) that the velocity $v[p]$ exists (for all $p$).  The monotonicity claim now follows from Theorem \ref{thm:other_mono}.\Qed

Note that one should be able to achieve {\em strict} monotonicity of the speed for the above model in high dimensions by using the method of \cite{HSun10}.

\subsection{Coupling with 1-d multi-excited random walks.}
\label{sec:excited_couple}
For RWRE satisfying (\ref{standinghypothesis}) such that there is sometimes a drift in direction $u$ but never a drift in direction $-u$ the walk should be almost surely transient in direction $u$ with positive speed.  
The transience result can be proved by considering the accumulated drift in direction $u$ (which is non-decreasing in this case) and adapting arguments appearing for example in Zerner \cite{Zern06}, while we expect that a proof of the speed result requires the extension of more technical machinery (such as Kalikow's condition, or methods used for studying excited random walk) to our non-elliptic setting.  The authors have some such proofs in preparation. 

If we can live with additional strong assumptions, such as in the following Lemma (whose proof is left as an exercise), it is relatively simple to show that the speed is positive. 
\begin{LEM}
\label{lem:liminfspeed}
If $\exists\epsilon>0$ such that $\sum_{e\in\mc{E}} e\cdot u \omega_o(e)\ge\epsilon$ $\mu$-a.s. then $\lim_{n\ra \infty}n^{-1}X_n\cdot u>0$.
\end{LEM}
In this paper we make what one might call intermediate assumptions, essentially that with sufficiently large probability we have a sufficiently large local drift. This enables us to give a relatively simple coupling proof that the speed is positive, by appealing to results from \cite{HS_comb, Zern05, BS08}.  For convenience, we state the result for the direction $u=e_1+\dots+e_d$, and we denote by $\overline{\omega}_x\in \Z^d$ the vector with entries $\overline{\omega}_x^{[i]}=\omega_x(e_i)-\omega_x(-e_i)$. Therefore the quenched drift at $x$ in the direction $u$ is $\sum_{e\in\mc{E}}e\cdot u\,\omega_x(e)=\overline{\omega}_x^{[i]}\cdot u$.

\begin{THM}
\label{thm:excouple}
Fix $d\ge 2$ and let $u=e_1+\dots+e_d=(1,\dots,1)$.  Suppose that $\mu(\overline{\omega}_x\cdot u \ge 0)=1$ and $\mu(\overline{\omega}_x\cdot u \ge a)=b$ for some $a,b>0$.  If  $\frac{ab}{1-b}>1$ (resp.~$>2$) the RWRE $X$ is transient (resp.~has positive speed) in direction $u$, $P$-almost surely.
\end{THM}
\proof If $b=1$ the walk satisfies $\liminf_{n\ra \infty} n^{-1}X_n\cdot u>0$ almost surely. So suppose that $b\in (0,1)$. Then the conditional measures $\mu^+(\omega_o\in \bullet)=\mu(\omega_o\in \bullet\mid\overline{\omega}_x\cdot u \ge a)$ and $\mu^-(\omega_o\in \bullet)=\mu(\omega_o\in \bullet\mid\overline{\omega}_x\cdot u < a)$ are well defined.  Let $\{W^+_{j,k}\}_{j \in \Z, k \in \N}$ and $\{W^-_{j,k}\}_{j \in \Z, k \in \N}$ be independent random variables with laws $\mu^+$ and $\mu^-$ respectively.

For $j\in \Z$ let $B_j=\{x\in \Z^d:x\cdot u=j\}$, and let $\{G_{j,k}\}_{j\in \Z,k\in \mathbb{N}}$ and $\{U_{j,k}\}_{j\in \Z,k\in \mathbb{N}}$ be independent random variables with laws~$\sim \text{Geometric}(1-b)$ and $\sim U[0,1]$ respectively.  The random variables $G_{j,k}$ will indicate the numbers of previously unvisited vertices in $B_j$ we have to visit before finding the next new site such that $\overline{\omega}\cdot u < a$.  Let $(u_1,\dots,u_{2d})=(e_1,\dots,e_d,-e_1,\dots,-e_d)$. 

Let $X_0=o$.  Define $\omega_{o}=W^+_{o,1}$ if $G_{o,1}>1$ and $\omega_{o}=W^-_{o,1}$ otherwise.  
Given $\{X_j, \omega_{X_j}\}_{j\le n}$, let $Y_n=\sum X_n^{[k]}=X_n \cdot u$, and let $L_n(j)=|\{i\le n:Y_i=j\}|$ be the local time of $Y$ at $j$ up to time $n$.  Set $L_n=L_n(Y_n)$.  Then define
\[X_{n+1}-X_n=u_i, \quad \text{ if }\sum_{j=1}^{i-1}\omega_{X_n}(u_j)<U_{Y_n,L_n}\le \sum_{j=1}^i \omega_{X_n}(u_j), \quad i=1,\dots, 2d.\]
Given $\{X_k\}_{k\le n+1}$ and $\{\omega_{X_k}\}_{k\le n}$ let 
\[ \omega_{X_{n+1}}=\begin{cases}
\omega_{X_l}, &\text{ if }X_{n+1}=X_l, l<n+1\\
W^-_{Y_{n+1},r}, &\text{ if }X_{n+1}\notin\{X_0,\dots,X_{n}\}, \text{ and}\\
  &\{X_k\}_{k\le n+1} \text{ has visited $r$ distinct sites in $B_{Y_{n+1}}$, and $r\in\{ \sum_{u=1}^sG_{Y_{n+1},u}\}_{s\in \N}$}\\
W^+_{Y_{n+1},r}, &\text{ if }X_{n+1}\notin\{X_0,\dots,X_{n}\}, \text{ and}\\
  &\{X_k\}_{k\le n+1} \text{ has visited $r$ distinct sites in $B_{Y_{n+1}}$, and $r\notin\{ \sum_{u=1}^sG_{Y_{n+1},u}\}_{s\in \N}$}.
\end{cases}\]
The reader can check that $X$ is then a random walk in a random environment $\omega$ that is i.i.d.~with $\omega_o\sim \mu$.

Note that the increments of $Y$ are in $\{-1,1\}$ and that $\mP(Y_{n+1}-Y_n=1|\omega_{X_n})=\sum_{j=1}^d \omega_{X_n}(u_i)$.  For at least the first $G_{j,1}-1$ visits of $X$ to $B_j$, the environment seen by the walker has law $\mu^+$ (not necessarily independently, as the same site could be visited more than once).  Thus for at least the first $G_{j,1}-1$ visits of 
$Y$ to $j$, the next increment of $Y$ has probability at least $\sum_{i=1}^d\omega(u_i)=(1+\overline{\omega}_x\cdot u)/2\ge (1+a)/2$ of being $+1$.  On subsequent visits, independent of the history it has probability at least $\hlf$ of being +1.  

Now consider a random walk $Z$ on $\Z$, with $Z_0=0$, that evolves as follows.  Given that $Z_n=j$ and $|\{k\le n:Z_k=j\}|=r$, 
\eq
\lbeq{Zdefn}
Z_{n+1}-Z_n=\begin{cases}
1, & \text{ if }r<G_{j,1} \text{ and }U_{j,r}\le (1+a)/2\\
-1, & \text{ if }r<G_{j,1}  \text{ and } (1+a)/2< U_{j,r}\\
1, & \text{ if }r\ge G_{j,1}  \text{ and } U_{j,r}\le\frac{1}{2}\\
-1, & \text{ if }r\ge G_{j,1}  \text{ and } \frac{1}{2}< U_{j,r}
\end{cases}
\en
The reader can check that this has coupled $Z$ and $Y$ together so that for all $j,r$
if $Y$ goes left on its $r$th departure from $j$ then so does $Z$ (if $Z$ visits $j$ at least $r$ times).

The random walk $Z$ defined in \refeq{Zdefn} is a multi-excited random walk in a random cookie environment $\upsilon$ such that $\{\upsilon(i,\cdot)\}_{i\in \Z}$ are i.i.d. with $\upsilon(0,r)=(1+a)/2$ for $r<G_{0,1}$ and $\upsilon(0,r)=\frac{1}{2}$ for $r\ge G_{0,1}$ (i.e.~a Geometric$(1-b)$ number of cookies at each site).  By \cite{Zern05,BS08}, $Z$ is transient to the right (resp.~has positive speed) if and only if $\alpha=E[\sum_{k\ge 1}(2\upsilon(o,k)-1)]>1$ (resp.~$>2$).  Now
\[E\left[\sum_{i\ge 1}(2\upsilon(0,i)-1)\right]\le E\left[\sum_{i=1}^{G_{0,1}}a\right]=\frac{ab}{1-b},\]
so $Z$ satisfies the claim of the proposition.  By \cite[Theorem 1.3]{HS_comb}, under this coupling, if $Z$ is transient to the right then so is $Y$, moreover $\limsup_{n\ra \infty} \frac{Z_n}{n}\le \limsup_{n\ra \infty} \frac{Y_n}{n}$ and the result follows (using Theorem \ref{thm:dtransv} to conclude that $Y$ has a limiting speed).
\Qed

In particular, from Theorem \ref{thm:excouple} we conclude that the uniform RWRE has positive speed in direction $u=(1,1)$ in the following cases:  $(\NSE,\leftrightarrow)$ and $(\NSE$, $\NSEWalt$ $)$ when $p>\frac{6}{7}$ ($a=\frac13$ and $b=p$); and $(\rightarrow$, $\NSEWalt$ $)$ and $(\NE$, $\NSEWalt$ $)$ when $p>\frac{2}{3}$ ($a=1$ and $b=p$).

\section{Calculation of speeds for uniform RWRE}
\label{sec:renewal}
There are many RWRE models for which it is actually obvious that transience holds and that speeds exist, due to the presence of a simple renewal structure (via a forbidden direction). In a number of cases, speeds can be calculated explicitly, and continuity and other properties of this speed as a function of certain parameters can be observed.  Example \ref{exa:3valuedspeed} is one such case, but this is somewhat messy.  The calculations all make use of the following lemma, whose proof is a simple application of the LLN/Renewal Theorem.
\begin{LEM} 
\label{lem:trivialtransience}
Assume (\ref{standinghypothesis}) and suppose that $\mu(\downarrow\in \mc{G}_o)=0$ but $\mu(\uparrow\in \mc{G}_o)>0$. 
Then the RWRE is transient in direction $e_2$, $P$-almost surely. 
Let $T$ be the first time the RWRE follows direction $e_2$. If $E[T]<\infty$ then $X_n$ has an asymptotic speed 
$v=(v^{[1]}, \dots,v^{[d]})$, in the sense that $P(n^{-1}X_n\to v)=1$. Moreover, $v^{[i]}=E[X^{[i]}_T]/E[T].$
\end{LEM}
The class of 2-valued uniform RWRE in 2-dimensions with a forbidden direction is a collection of relatively simple examples where $v$ can be calculated and continuity (as a function of $p$) and slowdown observed.  We will sketch the argument in the case of $(\NE,\NW)$ and will give a table, summarizing the results we know of in other 2-valued 2-dimensional models. Readers are referred to \cite{HS_speeds} for the detailed calculations in other cases.  



\begin{LEM}
\label{lem:NE_NWspeed}
Consider the uniform RWRE model ($\NE\NW$), i.e.~$\mu(\{\uparrow,\ra\})=p$ and $\mu(\{\uparrow, \la\})=1-p$. The asymptotic speed is $(v^{[1]}, v^{[2]})$ with 
\[v^{[1]}=\frac{(2p-1)(p^2-p+6)}{6(2-p)(1+p)},
\qquad v^{[2]}=\frac12.\]
\end{LEM}
\proof
For $n\ge 0$, let $\tau_n=\inf\{m\ge 0:X_m^{[2]}=n\}$.  Then  for $i\ge 1$, $T_i=\tau_i-\tau_{i-1}$ are i.i.d.~Geometric$(1/2)$ random variables (with mean $2$), and $Y_i=X^{[1]}_{\tau_i -1}-X^{[1]}_{\tau_{i-1}}$ are i.i.d.~random variables, independent of the $\{T_i\}_{i\ge 1}$.  So $E[T_i]=2$ and $v^{[2]}=1/2$. 
As in Lemma \ref{lem:trivialtransience} we have (almost surely as $n\ra \infty$)
\[\frac{Y_n^{[1]}}{n}\ra \frac{E[Y_1]}{E[T_1]}=\frac{E[Y_1]}{2}.\]
Letting $Y=Y_1$, it remains to calculate $E[Y]$. 

For $j\ge 1$, we can have $Y=j$ three ways -- reaching no $\NW$ vertex, reaching a $\NW$ vertex at $(j,0)$, or reaching a $\NW$ vertex at $(j+1,0)$. Thus 
\begin{align*}
P(Y=j)&=p^{j+1}\big(\frac12\big)^{j+1}+p^j(1-p)\sum_{n=0}^\infty\big(\frac12\big)^{j+2n+1}+p^{j+1}(1-p)\sum_{n=0}^\infty\big(\frac12\big)^{j+2n+3}\\
&=\frac{p^j(4-p^2)}{3\cdot 2^{j+1}}.
\end{align*}
Likewise, we can have $Y=-j$, $j\ge 1$ three ways, depending on where if anywhere $X_n$ reaches a $\NE$ vertex, giving 
$P(Y=-j)=\big((1-p)^j(4-(1-p)^2)\big)/\big(3\cdot 2^{j+1}\big).$ The case $j=0$ would be similar, but is not needed. Summing over $j$ gives that
\begin{align*}
E[Y]&=\frac{p(4-p^2)}{12}\cdot \frac{1}{(1-p/2)^2}-\frac{(1-p)(4-(1-p)^2)}{12}\cdot\frac{1}{(1-(1-p)/2)^2}\\
&=\frac{p(2+p)}{3(2-p)}-\frac{(1-p)(3-p)}{3(1+p)}
=\frac{(2p-1)(p^2-p+6)}{3(2-p)(1+p)}.
\end{align*}
\Qed

\begin{table}
\begin{center}
\begin{tabular}{l|l|l}
$\gamma^{\sss 1}, \gamma^{\sss 2}$ & Random walk  & Reference  \\
\hline
$\uparrow$ $\rightarrow$ & $v=(1-p,p)$.  &  As in Lemma \ref{lem:NE_NWspeed} \\
$\uparrow$ $\downarrow$ & Stuck on two vertices.   &  Lemma \ref{lem:criteriontobeinfinite} \\
$\leftrightarrow$ $\uparrow$ &  $v=\Big(0,\frac{(1-p)^2}{p+(1-p)^2}\Big)$. &  As in Lemma \ref{lem:NE_NWspeed} \\
$\leftrightarrow$ $\rightarrow$ &  $v=\Big(\frac{1-p}{1+p}, 0\Big)$.  &  As in Lemma \ref{lem:NE_NWspeed}\\
$\leftrightarrow$ $\updownarrow$ &  $v=(0,0)$. &   $\text{Symmetry}^1$\\
\hline
$\NE$ $\uparrow$ &  $v=\Big(\frac{p}{2}, 1-\frac{p}{2}\Big)$. &  As in Lemma \ref{lem:NE_NWspeed}\\
$\NE$ $\NW$ &  $v=\Big(\frac{(2p-1)(p^2-p+6)}{6(2-p)(1+p)},\hlf\Big)$.  &   As in Lemma \ref{lem:NE_NWspeed}\\
$\NE$ $\leftrightarrow$ &  $v=\left(\frac{1}{p^2}+\frac{(1-p)^2}{2p(1-p+p\log p)}\right)^{-1}\cdot (1,1)$.  & As in Lemma \ref{lem:NE_NWspeed}\\
$\NE$ $\leftarrow$ &  $v=\left(\frac{p(2-p)}{2+3p-2p^2-p^3}\right)\cdot (3,1)+(-1,0)$.  &  As in Lemma \ref{lem:NE_NWspeed} \\
$\NE$ $\SW$ & $v^{[1]}=v^{[2]}\uparrow$ in $p$. Transient${}^2$ for $p\approx 0,1$.   &   Cor. \ref{cor:orthant_D_trans} / Thm. \ref{thm:other_mono}
\\
&{\bf Conjecture:} $v\ne 0$ for $p\neq\hlf$, Recurrent when $p=\hlf$  &\vspace{.1cm} \\
\hline
$\SWE$ $\downarrow$ 
& 
$\frac{1}{v^{[2]}}=\frac{8p(1-p)}{1+\sqrt{5}}-1-2p-\frac{4(1-p)^2(5+\sqrt{5})}{p(1+\sqrt{5})}\overset{\infty}{\underset{n=2}{\sum}}\frac{p^k}{1+2^{-k}(3+\sqrt{5})^k}$.  & As in Lemma \ref{lem:NE_NWspeed}${}^4$\\
$\SWE$ $\rightarrow$ & $-\frac{1}{v^{[2]}}=4-p-\frac{5+\sqrt{5}}{2}(1-p)^2\Theta(p\gamma)
+\frac{(1-p)[3+\sqrt{5}-(1-p)(5+\sqrt{5})\Theta(p)]^2}{(3+\sqrt{5})[2-(1-p)(5+\sqrt{5})\Theta(p\gamma)]}$ & As in  Lemma \ref{lem:NE_NWspeed}${}^4$\\
& where $\gamma=\frac{3+\sqrt{5}}{2}$ and $\Theta(z)=\sum_{n=0}^\infty \frac{z^n}{\gamma^{2n+1}+1}$. $v^{[1]}=1-3v^{[2]}$. 
 & \\
$\SWE$ $\uparrow$ & $v^{[1]}=0$, $v^{[2]}\downarrow$ in $p$. Transient${}^2$ for $p\approx 0$.  &   Cor. \ref{cor:Mdirect}\\
& {\bf Conjecture: } $\exists ! p(\ne 3/4)$ s.t. $v[p]=0$. Recurrent for this $p$. & \\
$\SWE$ $\leftrightarrow$ & $v^{[1]}=0$, $v^{[2]}<0$ for $p>0$. $v^{[2]}$ strictly $\downarrow$ in $p$. & As in  Lemma \ref{lem:NE_NWspeed}${}^5$\\
$\SWE$ $\updownarrow$ & $v^{[1]}=0$, $v^{[2]}\downarrow$ in $p$. Transient${}^3$ for $p>\frac{3}{4}$,  $v^{[2]}< 0$ for $p>\frac{6}{7}$. & Thm. \ref{thm:other_mono} / Thm.  \ref{thm:excouple}\\
&{\bf Conjecture:} $v^{[2]}< 0$ for $p>0$.   & \\
$\SWE$ $\NE$ & $3v^{[2]}=5v^{[1]}-1$. $v^{[1]}\downarrow$ in $p$. & Thm. \ref{thm:other_mono} / Cor. \ref{cor:vrelation}\\
$\SWE$ $\SW$ & $v^{[1]}=1+3v^{[2]}$   & As in Lemma \ref{lem:NE_NWspeed}${}^6$\\
$\SWE$ $\NSE$ & $v\cdot(1,-1)=\frac{1}{3}$, $v\cdot (1,1)\downarrow$ in $p$.  & Thm. \ref{thm:other_mono} / Cor. \ref{cor:vrelation}\\
$\SWE$ $\NWE$ &  $v^{[1]}=0$, $v^{[2]}\downarrow$ in $p$. & Thm. \ref{thm:other_mono} / Cor. \ref{cor:vrelation} \\
& {\bf Conjecture:} $v^{[2]}\ne 0$ for $p\neq\frac12$. Recurrent when $p=\frac12$.&   \vspace{.1cm}\\
\hline 
$\NSEWalt$ \vphantom{${A^A}^A$} \hspace{.1cm} $\uparrow$ & $v^{[1]}=0$, $v^{[2]}\downarrow$ in $p$. Transient${}^3$ for $p<\hlf$, $v^{[2]}> 0$ for $p<\frac13$.  & Thm. \ref{thm:other_mono} / Thm.  \ref{thm:excouple}\\
&  {\bf Conjecture:} $v^{[2]}> 0$ for $p<1$. & \\
$\NSEWalt$ \hspace{.1cm} $\NE$ & $v^{[1]}=v^{[2]}\downarrow$ in $p$. Transient${}^3$ for $p<\hlf$, $v^{[1]}> 0$ for $p<\frac13$.  & Thm. \ref{thm:other_mono} / Thm.  \ref{thm:excouple}\\
&  {\bf Conjecture:}  $v^{[1]}> 0$ for $p<1$.  & \\
$\NSEWalt$ \hspace{.1cm} $\WE$ &  $v=(0,0)$  & Symmetry${}^1$.  \\
$\NSEWalt$ \hspace{.1cm} $\SWE$ & $v^{[1]}=0$, $v^{[2]}\uparrow$ in $p$. 
Transient${}^3$ for $p<\frac{1}{4}$,  $v^{[2]}< 0$ for $p<\frac{1}{7}$. & Thm. \ref{thm:other_mono} / Thm.  \ref{thm:excouple}\\
&{\bf Conjecture:} $v^{[2]}< 0$ for $p<1$.&\\
\hline
\end{tabular}
\end{center}
\caption{Table of results for uniform RWRE in 2-dimensional 2-valued degenerate random environments, where the first configuration occurs with probability $p\in (0,1)$ and the other with probability $1-p$.}
\label{tab:walks}
\end{table}

\medskip

Table \ref{tab:walks} summarizes what we know about uniform RWRE in 2-dimensional 2-valued random environments. Explicit speeds are calculated as in Lemma \ref{lem:NE_NWspeed} (for details, consult \cite{HS_speeds}). All other conclusions follow immediately from results stated in the paper. Note that many of the conjectures would follow if we knew that speeds were continuous in $p$ and that monotonicity was strict.

Note that there is a related table in \cite{HS_DRE1}, giving percolation properties for the directed graphs $\mc{C}$ and $\mc{M}$. The latter includes 2-valued environments such as $(\NSEWalt \,\,\, , \cdot)$ (site percolation), in which one of the possible environments has no arrows. These environments do not appear in the present table, because (as remarked in Section \ref{sec:model}), the walk gets stuck on a finite set of vertices (in this case 1 vertex) the RWRE setup we have chosen requires that motion be possible in at least one direction. \medskip

\noindent {\bf Notes to Table \ref{tab:walks}} \newline
${}^1$ As indicated following Lemma \ref{lem:Mrecurrent}, it follows from results of Berger \& Deuschel \cite{BD11} that $\mc{M}$ is recurrent $\forall p$. \newline
${}^2$ Bounds on the
critical probability are given in \cite{HS_DRE1}. Improved bounds are in preparation.\newline
${}^3$ Improved ranges of values giving transience and speeds are in preparation. \newline
${}^4$ Martin Muldoon has pointed out that this can also be expressed in terms of $q$-hypergeometric functions.\newline
${}^5$ We do not have a closed form expression for this. But asymptotic expressions are feasible.\newline
${}^6$ An expansion as in the case ($\SWE$ $\rightarrow$) should be possible.

\section*{Acknowledgements}
Holmes's research is supported in part by the Marsden fund, administered by RSNZ. Salisbury's research is supported in part by NSERC. Both authors acknowledge the hospitality of the Fields Institute, where part of this research was conducted.


\begin{thebibliography}{10}

\bibitem{BS08}
A.-L.~Basdevant and A.~Singh.
\newblock On the speed of a cookie random walk.
\newblock {\em Probab. Theory Relat. Fields.}, 141(3-4):625--645, 2008.


\bibitem{BD11}
N.~Berger and J.-D.~Deuschel.
\newblock A quenched invariance principle for non-elliptic random walk in I.I.D. balanced random environment.
\newblock Preprint, 2011.


\bibitem{BSZ03}
E.~Bolthausen, A.-S. Sznitman, and O.~Zeitouni.
\newblock Cut points and diffusive random walks in random environment.
\newblock {\em Ann. Inst. H. Poincar\'e Probab. Statist.}, 39(3):527--555,
  2003.


\bibitem{HS_DRE1}
M.~Holmes and T.S.~Salisbury.
\newblock Degenerate random environments.
\newblock Preprint, 2009.

\bibitem{HS_DRE2}
M.~Holmes and T.S.~Salisbury.
\newblock Degenerate random environments II.
\newblock Preprint, 2009.

\bibitem{HS_comb}
M.~Holmes and T.S.~Salisbury.
\newblock A combinatorial result with applications to self-interacting random
  walks.
\newblock J. Combinatorial Theory A, 19: 460--475 (2012).

\bibitem{HS_speeds}
M.~Homes and T.S.~Salisbury.
\newblock Speed calculations for random walks in degenerate random environments. 
\newblock Unpublished notes (2011). {\tt www.math.yorku.ca/~salt/preprints/speeds2011May31.pdf}

\bibitem{HSun10}
M.~Holmes and R.~Sun.
\newblock A monotonicity property for random walk in a partially random environment.
\newblock Preprint, 2010.

\bibitem{K81}
S.A.~Kalikow. 
\newblock Generalized random walks in random environment. 
\newblock {\em Ann. Probab.}, 9: 753--768, 1981.

\bibitem{Law82}
G.~F. Lawler.
\newblock Weak convergence of a random walk in a random environment.
\newblock {\em Comm. Math. Phys.}, 87(1):81--87, 1982/83.

\bibitem{MadrasT}
N.~Madras and D.~Tanny.
\newblock Oscillating random walk with a moving boundary.
\newblock{\em Israel J. Math.}, 88:333-365, 1994.

\bibitem{RS06}
F.~Rassoul-Agha and T.~Sepp\"{a}l\"{a}inen.
\newblock Ballistic random walk in a random environment with a forbidden direction.
\newblock {\em ALEA}, 1:111--147, 2006.


\bibitem{SZ99}
A.-S. Sznitman and M.~Zerner.
\newblock A law of large numbers for random walks in random environment.
\newblock {\em Ann. Probab.}, 27:1851--1869, 1999.

\bibitem{Zeit04}
O.~Zeitouni.
\newblock Random walks in random environment.
\newblock In {\em Ecole d'Et\'e de Probabilit\'es de Saint Flour 2001, Lecture
  Notes in Mathematics, no. 1837}. Springer-Verlag, Berlin, 2004.

\bibitem{Zern02}
M.P.W.~Zerner.
\newblock A non-ballistic law of large numbers for random walks in i.i.d. random environment.
\newblock {\em Electron. Comm. Probab.}, 7:191--197 (electronic), 2002.

\bibitem{Zern05}
M.P.W.~Zerner.
\newblock Multi-excited random walks on integers.
\newblock {\em Probab. Theory Relat. Fields.}, 133:98--122, 2005.

\bibitem{Zern06}
M.P.W.~Zerner.
\newblock Recurrence and transience of excited random walks on {$\Z\sp d$}
  and strips.
\newblock {\em Electron. Comm. Probab.}, 11:118--128 (electronic), 2006.

\bibitem{Zern07}
M.P.W.~Zerner.
\newblock The zero-one law for planar random walks in i.i.d. random environments revisited. 
\newblock {\em Electron. Comm. Probab.}  12:326?-335 (electronic), 2007.


\bibitem{MZ01}
M.P.W.~Zerner and  F.~Merkl.
\newblock A zero-one law for planar random walks in random environment. 
\newblock {\em Ann. Probab.} 29(4):1716--1732, 2001. 


\end{thebibliography}
\bibliographystyle{plain}

\end{document}